\documentstyle[a4wide]{article}
\begin{document}
\newcommand{\qed}{\hfill $\Box$ \\[0.2cm]}
\newcommand{\dkz}{{\bf Proof:}{\hspace{0.3cm}} \nopagebreak}
\newcommand{\HDS}{\vrule width0pt height2.3ex depth1.05ex\displaystyle}
\newcommand{\cls}{${\cal C}(\pi_1,\pi_2)$}
\newcommand{\G}{$\cal G$}
\newcommand{\Gr}{${\cal G}^r$}
\newcommand{\Gz}{${\cal G}^z$}
\newcommand{\str}{\rightarrow}
\newcommand{\I}{\wedge}
\def\teo#1{{\bf Theorem #1}\hspace{1em}}
\def\lema#1{{\bf Lemma #1}\hspace{1em}}
\def\Pime#1{\makebox[3ex][l]{\mbox{\em #1}}}
\def\tekst#1{\mbox{\rm #1}}
\def\f#1#2{{{\HDS #1}\over{\HDS #2}}}
\def\fp#1#2{{{\HDS #1}\atop{\HDS #2}}}
\def\afrac#1{{\phantom{\HDS #1}\atop{\HDS #1}}}
\def\bfrac#1#2{{\phantom{\HDS \Theta}\atop{\phantom{
\vrule width0pt height2.3ex depth#2\displaystyle \Theta}\atop{\HDS #1}}}}
\def\rza{{\mbox{\hspace{0.5em}}}}
\def\rzb{{\mbox{\hspace{2em}}}}
\def\rzbc{{\mbox{\hspace{4,5em}}}}
\def\rzc{{\mbox{\hspace{7em}}}}
\def\rzd{{\mbox{\hspace{11em}}}}
\def\Apak#1#2#3#4{\f{\f{\fp{#1}{#2}}{\cdots} \pravila{#3}}{#4}} 
\def\Bpak#1#2#3{\f{\f{#1}{\cdots} \pravila{#2}}{#3}} 
\def\Cpak#1#2#3{\f{\f{#1}{\cdots} \ppravila{#2}}{#3}} 
\def\Dpak#1#2#3{\f{\f{#1}{\cdots} \pppravila{#2}}{#3}} 
\def\pravilo#1{ \rza \makebox[-.5em][l]{\mbox{\rm #1}}}
\def\lpravilo#1{ \makebox[-.5em][r]{\mbox{\rm #1}} \rza}
\def\pravila#1{ \rza \makebox[-.5em][l]{\raisebox{-1.7ex}
{\parbox{15ex}{\footnotesize \baselineskip=0.5\baselineskip #1}}}}
\def\ppravila#1{ \rza \makebox[-.5em][l]{\raisebox{-1.7ex}
{\parbox{19ex}{\footnotesize \baselineskip=0.5\baselineskip #1}}}}
\def\pppravila#1{ \rza \makebox[-.5em][l]{\raisebox{-1.7ex}
{\parbox{30ex}{\footnotesize \baselineskip=0.5\baselineskip #1}}}}
\def\naz#1{$\tekst{#1}^z$}
\title
{On Permuting Cut with Contraction}

\author{ }

\date{ }
\maketitle
\begin{abstract}
This paper presents a cut-elimination procedure for intuitionistic
propositional logic in which cut is eliminated directly, without
introducing the multiple-cut rule mix, and in which pushing cut above
contraction is one of the reduction steps. The presentation of this
procedure is preceded by an analysis of Gentzen's mix-elimination
procedure, made in the perspective of permuting cut with contraction.
It is also shown that in the absence of implication, pushing cut above
contraction doesn't pose problems for directly eliminating cut.
\end{abstract}
\vspace{0.5cm}
\section{Introduction}
\noindent The structural rule of contraction poses special problems
for cut elimination.  It is because of contraction that in the
cut-elimination procedure of [1935] Gentzen replaced his rule
\[ \f{\Gamma \vdash \Theta, A \rzb A, \Delta \vdash \Lambda}
{\Gamma, \Delta \vdash \Theta, \Lambda} \pravilo{(Gentzen's cut)} \]
by a rule derived from cut, contraction and interchange, called {\em mix}
({\em Mischung} in German),
\[ \f{\Gamma \vdash \Theta \rzb  \Delta \vdash \Lambda}
{\Gamma, \Delta^{\ast} \vdash \Theta^{\ast}, \Lambda}  \]
where $\Theta$ and $\Delta$ are sequences of formulae in each of
which occurs at least once a formula $A$, called the {\em mix-formula},
and $\Theta^{\ast}$ and $\Delta^{\ast}$ are obtained from, respectively,
$\Theta$ and $\Delta$ by deleting all occurrences of $A$. That cut can be
eliminated is then demonstrated by eliminating mix.

Mix also solves a problem involving the structural rule of interchange.
Namely, we cannot permute (Gentzen's cut) with an interchange above the cut
involving the cut formula $A$, because Gentzen required that the cut formula
be the last formula in the sequence on the right-hand side of the left
premise and the first formula on the left-hand side of the right premise.
However, this problem is easily solved by replacing (Gentzen's cut) with
\[ \f{\Gamma \vdash \Theta_{1}, A, \Theta_{2} \rzb \Delta_{1}, A,
\Delta_{2} \vdash \Lambda}
{\Delta_{1}, \Gamma, \Delta_{2} \vdash \Theta_{1}, \Lambda, \Theta_{2}}
\pravilo{(cut)} \]
which in the context of classical and intuitionistic logic
doesn't represent an essential departure from the original 
systems of Gentzen. We can always obtain the effect of 
the rule (cut) with the help
of (Gentzen's cut) and interchanges preceding and following this cut.

The special problem brought for cut elimination by contraction, because
of which mix is introduced, occurs when we have to permute a cut with
contraction above the cut involving the cut formula $A$, i.e. when we
push a cut above such a contraction. If a figure with a topmost cut
$$ \f
{
\afrac{\Gamma \vdash \Theta, A}
\rzb
\f{A, A, \Delta \vdash \Lambda}{A, \Delta \vdash \Lambda} \pravilo{contraction}
}
{\Gamma, \Delta \vdash \Theta, \Lambda} \pravilo{cut}
\leqno(\ast)
$$
is replaced by the figure
$$ \f
{
\afrac{\Gamma \vdash \Theta, A}
\rzb
\f
{\Gamma \vdash \Theta, A \rzb A, A, \Delta \vdash \Lambda}
{\Gamma, A, \Delta \vdash \Theta, \Lambda} \pravilo{cut}
}
{
\Bpak{\Gamma, \Gamma, \Delta \vdash \Theta, \Theta, \Lambda}
{interchanges and contractions}{\Gamma, \Delta \vdash \Theta, \Lambda}
}
\pravilo{cut}
\leqno(\ast\ast)
$$
we have two cuts with the same cut formula $A$ replacing a single
cut with this cut formula. Of these two cuts, the upper cut has lower
rank and can be eliminated by the induction hypothesis, but after
this elimination is made, the remaining, lower, cut, which has now
become topmost, need not have lower rank than the original cut.

On the other hand, if the following figure with a topmost mix
$$ \f
{
\afrac{\Gamma \vdash \Theta}
\rzb
\f{A, A, \Delta \vdash \Lambda}{A, \Delta \vdash \Lambda} \pravilo{contraction}
}
{\Gamma, \Delta^{\ast} \vdash \Theta^{\ast}, \Lambda} \pravilo{mix}
$$
is replaced by the figure
\[ \f{\Gamma \vdash \Theta \rzb  A, A, \Delta \vdash \Lambda}
{\Gamma, \Delta^{\ast} \vdash \Theta^{\ast}, \Lambda} \pravilo{mix}  \]
then the new, single, application of mix is topmost and has lower rank
than the original mix.

It is sometimes assumed that Gentzen's cut-elimination procedure
is based on replacing $(\ast)$ by $(\ast\ast)$ (see [Carbone 1997], p. 285).
When mix is reconstructed in terms of cut and other structural rules,
Gentzen's procedure does indeed involve pushing cut above contraction,
but only as part of more complicated steps, as we shall show in Section
2 below. It might even be taken that in some of these steps cut
is pushed below contraction, in the opposite direction.

When some forty years ago Lambek undertook in [1958] to eliminate cut
in a contractionless sequent system, he didn't need to bother with mix,
and could eliminate cut directly. Of course, one also need not rely on mix
in other contractionless systems of substructural logics that have been
introduced since: namely, systems of BCK logic and linear logic.

In [1978], Szabo attempted to systematize the cut-elimination algorithm
so that it can apply to a number of systems, with and without contraction.
In this algorithm, when contraction is present, cut is permuted with
contraction by passing from a figure like $(\ast)$ to a figure like
$(\ast\ast)$ (see [Szabo 1978], Appendix C, C.19.3, p. 234, C.38.3,
p. 239). To demonstrate that the figure of $(\ast\ast)$ is somehow simpler,
Szabo introduced in [1978] (pp. 242-243) a measure of complexity
counting the number of contractions above a cut. However, Szabo's measure
fails to show that the lower cut in $(\ast\ast)$ will have a smaller
measure of complexity, as can be seen in a counterexample presented in detail
in the last section of [B. 1997].\footnote{We are grateful to Andreja
Prijatelj for pointing a long time ago to one of us that Szabo's treatment of
the matter is unsatisfactory.}

Actually, one cannot push a cut above both a contraction on the left and a
contraction on the right, as the following simple counterexample shows. The
figure
$$ \f
{
\lpravilo{contraction}\f{\Gamma \vdash \Theta, A, A}{\Gamma \vdash \Theta, A}
\rzb
\f{A, A, \Delta \vdash \Lambda}{A, \Delta \vdash \Lambda} \pravilo{contraction}
}
{\Gamma, \Delta \vdash \Theta, \Lambda} \pravilo{cut}
$$
is not replaceable by a figure where all the cuts will be above all
contractions. Szabo doesn't eschew problems posed by this figure, though
he requires in [1978] (p. 234) that the right rank of the cut be 1 if we
want to diminish the left rank. Because, after permuting the cut in the
figure with the contraction on the left above the right premise, we obtain
$$ \f
{
\lpravilo{contraction}
\f{\Gamma \vdash \Theta, A, A}{\Gamma \vdash \Theta, A}
\rzc
\f
{\lpravilo{contraction}\f{\Gamma \vdash \Theta, A, A}
{\Gamma \vdash \Theta, A} \rzb \afrac{A, A, \Delta \vdash \Lambda}}
{\Gamma, A, \Delta \vdash \Theta, \Lambda} \pravilo{cut}
}
{
\Bpak{\Gamma, \Gamma, \Delta \vdash \Theta, \Theta, \Lambda}
{interchanges and contractions}{\Gamma, \Delta \vdash \Theta, \Lambda}
}
\pravilo{cut}
$$
where the upper cut may be of the lowest possible right rank. When
we next permute this upper cut with contraction on the right over the left
premise, we obtain the figure 
$$ \f
{
\bfrac
{\lpravilo{contraction}\f{\Gamma \vdash \Theta, A, A}
{\Gamma \vdash \Theta, A}}{3.5ex}
%\bfrac{\jedan}
%{\bfrac{\jedan}{\lpravilo{contraction}\f{\Gamma \vdash \Theta, A, A}
%{\Gamma \vdash \Theta, A}}}
\rzb
\f
{\lpravilo{cut}\f{\Gamma \vdash \Theta, A, A \rzb A, A, \Delta \vdash \Lambda}
{\Gamma, A, \Delta \vdash \Theta, A, \Lambda} \rzb
\afrac{A, A, \Delta \vdash \Lambda}}
{\Bpak{\Gamma, A, \Delta, A, \Delta \vdash \Theta, \Lambda, \Lambda}
{interchanges and contractions}{\Gamma, A, \Delta \vdash \Theta, \Lambda}
}
\pravilo{cut}
}
{
\Bpak{\Gamma, \Gamma, \Delta \vdash \Theta, \Theta, \Lambda}
{interchanges and contractions}{\Gamma, \Delta \vdash \Theta, \Lambda}
}
\pravilo{cut}
$$
where the lowest cut is in the same position as the initial one.

However, this does not exclude that Szabo's complexity measure could be
replaced by another measure, presumably more complicated, which would
show that the algorithm he envisaged would terminate if we have only
contraction on the left in a system close to Gentzen's $LJ$ of [1935], i.e., in
a system for intuitionistic logic. (In [Girard et al. 1992] something
like this measure is computed, but the elimination of all cuts is not
sought: in particular, some difficult cuts with contracted cut formulae
are not eliminated.) These matters are very much tied to the particular
formulation of a system. Zucker shows in [1974] (\S 7) that if in a system
for intuitionistic logic one replaces Gentzen's ``additive'', i.e.
{\em lattice}, rules for disjunction by ``multiplicative'' rules,
a procedure such as envisaged by Szabo would not terminate.

Actually, it is not difficult to find such a measure in the absence of
implication, as we show below, in Section 4. The presence of implication
poses special problems, for which we shall devise a cut-elimination procedure
that involves permuting contractions with other rules, and not only
with cut. Such permutations of contraction were studied in [Kleene 1952],
[Zucker 1974], [Minc 1996] and [Dyckhoff \& Pinto 1997], but we are not
aware that they have been integrated before into a cut-elimination
procedure. (Among these papers only Zucker's envisages 
permuting contraction with cut.)

The goal of this paper is to present this cut-elimination procedure for
intuitionistic propositional logic, in which cut is directly eliminated,
without passing via mix, and in which pushing cut above contraction, i.e.
passing from $(\ast)$ to $(\ast\ast)$, is a 
reduction step. The cut-elimination
procedure of [B. 1997] also eliminates cut directly, and it involves
pushing cut above contraction, but it is different and more entangled
than the procedure
we are going to present here. In a procedure envisaged by [Carbone 1997],
reminiscent of Curry's mix-elimination procedure (see [Curry 1963],
Chapter 5, D2), cut should be directly eliminated, but without pushing it
above contraction.

Although we suppose our procedure could be extended to the whole of
intuitionistic predicate logic, we restrict ourselves to the propositional
case, to make the exposition simpler. Anyway, our result is rather of
theoretical, and not of practical, interest. If one is just interested in
eliminating cut, and does not care how exactly this is done, 
Gentzen's solution
based on mix is simpler. It is probably optimal.

However, our procedure may perhaps come in handy in studies of complexity
of proofs. It exhibits more clearly than Gentzen's procedure that
contraction is the culprit for the hyperexponential growth of proofs in
cut elimination.

Anyway, it seems worth knowing that cut can be eliminated by pushing it
above contraction. If for nothing else, then to block inept criticism that
would confuse ``I don't know how to eliminate cut by pushing it above
contraction'' with ``Cut cannot be so eliminated''.

Our procedure consists of three phases. In the first phase we push
contractions below all rules, including cut, except for the rule of
introduction of implication on the right. Proofs where this has been
accomplished are called ``W-normal''. In the second phase, to reduce
the rank, we push cuts above other rules, among which, because of
W-normality, we don't have any more troublesome applications of contraction,
like those in $(\ast)$. This phase involves essentially permuting cuts
with cuts, which is a matter only implicitly and incompletely present in
Gentzen's procedure (see the comments below ($2^{\ast}$) in Section 2,
and cases (2.4), (3.7) and (3.8) in the proof of Theorem 6.1; see also
the passage from ($3$P) to ($3^{\ast}$P) in Section 2, the end of
Section 2 and the beginning of Section  7). However, this permuting is
prominent in categorial proof theory: it corresponds to associativity
of composition and to bifunctoriality equalities. In the third phase,
we reduce cuts to cuts of lower degree. Then we reenter the first
phase of W-normalizing, and then again we go into the second phase, etc.
The last phase will be a second phase where only cuts with axioms
remain, and these are then eliminated.

Before describing precisely this procedure we consider in the next section
(Section 2) what Gentzen's mix-elimination procedure has to say about
permuting cut with contraction when the mix rule is reconstructed in terms
of cut, contraction and interchange. In Section 3 we introduce formally our
variant of Gentzen's sequent system $LJ$ of intuitionistic propositional
logic, which we call \G. The main difference between $LJ$ and \G\  is that
in the latter we have rules like (cut) above, instead of (Gentzen's cut).
In Section 4 we show by a simple argument that in implicationless \G\  we can
eliminate cut by freely pushing cuts above contractions. Perhaps, as
Szabo supposed, such a free policy of pushing cut above contraction leads
to cut elimination in \G\  even in the presence of implication, 
but we have been
unable to show that indeed it does.

In the last two sections we present our cut-elimination procedure. Section
5 is devoted to W-normalizing, and Section 6 to the remaining phases of the
procedure. In Section 7 we make some concluding comments.

\section{Cut elimination via mix elimination in $LJ$}

\noindent Gentzen's mix rule is derivable in the presence of the
structural rules of cut, contraction and interchange. However, for any mix
of Gentzen's system $LJ$ of [1935]
\[\f{\Gamma \vdash A \rzb \Delta \vdash \Lambda}
{\Gamma, \Delta^{\ast} \vdash \Lambda}\pravilo{(mix)} \]
where $A$ happens to occur in $\Delta$ more than once, there is no unique
way to reconstruct it in terms of cut, contraction and interchange. For
example, the following instance of (mix)
\[\f{B, C \vdash A \rzb A, A, D, A \vdash E}
{B, C, D \vdash E}\pravilo{mix} \]
can be reconstructed either as a number of cuts and interchanges followed
by contractions:
\[ \f{\bfrac{B, C \vdash A}{8.3ex} \rzb
\f{B, C \vdash A \rzb A, A, D, A \vdash E}{\Bpak{B, C, A, D, A \vdash E}
{interchanges}{A, B, C, D, A \vdash E}}\pravilo{$LJ$ cut}}
{\f{\afrac{B, C \vdash A} \rzb \Bpak{B, C, B, C, D, A \vdash E}
{interchanges}{A, B, C, B, C, D \vdash E}}
{\Bpak{B, C, B, C, B, C, D \vdash E}{interchanges and contractions}
{B, C, D \vdash E}}\pravilo{$LJ$ cut}}\pravilo{$LJ$ cut}\]
or as interchanges and contractions followed by a single cut:
\[\f{\bfrac{B, C \vdash A}{8.3ex} \rzb \f{A, A, D, A \vdash E}
{\Bpak{A, A, A, D \vdash E}{contractions}{A, D \vdash E}}
\pravilo{interchange}}{B, C, D \vdash E}\pravilo{$LJ$ cut}\]
or in many other ways intermediate between these two extremes, as for
instance
\[ \f{\bfrac{B, C \vdash A}{8.3ex} \rzb \f{\afrac{B, C \vdash A} \rzb
\f{A, A, D, A \vdash E}{A, D, A \vdash E}\pravilo{contraction}}
{\Bpak{B, C, D, A \vdash E}{interchanges}{A, B, C, D \vdash E}}
\pravilo{$LJ$ cut}}{\Bpak{B, C, B, C, D \vdash E}
{interchanges and contractions}{B, C, D \vdash E}}\pravilo{$LJ$ cut}\]
We call the first of these reconstructions, with many cuts,
{\em polytomic}, while the second, with a single cut, will be
{\em monotomic}. Note that in the polytomic reconstruction, and in the
intermediate third reconstruction, the left premise of mix appears more than
once. To pass from such reconstructions to the mix reconstructed, we have to
apply a contraction principle of higher level, which permits to omit
repetitions among the sequents that make the premises of a rule.

Note also that the polytomic reconstruction of a mix is not unique: one
such reconstruction may be obtained from another by 
introducing interchanges and by
permuting $LJ$ cuts with other $LJ$ cuts. The order of contractions
in the bottom of the reconstruction is also not uniquely determined. It is
possible to make this reconstruction unique by introducing an order among
the rules involved in the reconstruction, the shortest way being
to attack first the leftmost formula. However, there is something arbitrary
in this order.

Whether Gentzen's mixes of $LJ$ will be reconstructed polytomically,
monotomically or in some other, intermediate, way is a matter of choice.
This choice is of no consequence if the goal is just to eliminate cut
by whatever means. However, if we are interested in describing exactly the
cut-elimination procedure, and wish to reconstruct this procedure from the
mix-elimination procedure, we will not end up by the same algorithm
if we reconstruct mix always polytomically or always monotomically.

Let us now investigate when cut has to be pushed above contraction
involving the cut formula in the uniform polytomic and uniform
monotomic reconstructions; namely, in the reconstruction
where mixes are always reconstructed polytomically and
in the reconstruction where mixes are always reconstructed
monotomically. We shall only consider these uniform
reconstructions. (Note that passing from the monotomic to the
polytomic reconstruction of a mix may itself be conceived as obtained
by pushing cut above contraction.)

If the right rank of a mix is equal to 1, then this mix is just an $LJ$ cut.
So we have only to consider cases where the right rank of the mix is
greater than 1 (see [Gentzen 1935], Section III.3121). The first interesting
case for us is when we have
$$\f{\afrac{\Gamma \vdash A} \rzb \f{A, A, \Delta \vdash \Lambda}
{A, \Delta \vdash \Lambda}\pravilo{contraction}}
{\Gamma, \Delta^{\ast} \vdash \Lambda}\pravilo{mix}\leqno{(1)} $$
and there are $n$ occurrences of $A$ in $\Delta$. Polytomically, (1)
is reconstructed as
$$ \f{\bfrac{\Gamma \vdash A}{8.3ex} \rzb \f{\afrac{\Gamma \vdash A} \rzb
\f{A, A,\Delta \vdash \Lambda}{A, \Delta \vdash \Lambda}\pravilo{contraction}}
{\Cpak{\Gamma, \Delta \vdash \Lambda}{interchanges and $n-1$ applications of
$LJ$ cut}{A, \Gamma,\ldots,\Gamma, \Delta^{\ast} \vdash \Lambda}}
\pravilo{$LJ$ cut}}{\Bpak{\Gamma, \Gamma,\ldots,\Gamma, \Delta^{\ast}
\vdash \Lambda}{interchanges and contractions}{\Gamma, \Delta^{\ast}
\vdash \Lambda}}\pravilo{$LJ$ cut} \leqno{\raisebox{3ex}{(1P)}} $$
and monotomically as
$$\f{\bfrac{\Gamma \vdash A}{8.3ex} \rzb \f{A, A, \Delta \vdash \Lambda}
{\Bpak{A, \Delta \vdash \Lambda}{interchanges and contractions}
{A, \Delta^{\ast} \vdash \Lambda}}
\pravilo{contraction}}{\Gamma, \Delta^{\ast} \vdash \Lambda}
\pravilo{$LJ$ cut} \leqno{\raisebox{7ex}{(1M)}} $$
In [1935] (III.3.121.21) Gentzen transforms (1) into
$$ \f{\Gamma \vdash A \rzb A, A, \Delta \vdash \Lambda}
{\Gamma, \Delta^{\ast} \vdash \Lambda}\pravilo{mix} \leqno{(1^{\ast})} $$
Polytomically, ($1^\ast$) is reconstructed as
$$ \f{\bfrac{\Gamma \vdash A}{8.3ex} \rzb \f{\bfrac{\Gamma \vdash A}{8.3ex}
\rzb \f{\Gamma \vdash A \rzb A, A,\Delta \vdash \Lambda}
{\Bpak{\Gamma, A, \Delta \vdash \Lambda}{interchanges}
{A, \Gamma, \Delta \vdash \Lambda}}\pravilo{$LJ$ cut}}
{\Cpak{\Gamma, \Gamma, \Delta \vdash \Lambda}
{interchanges and $n-1$ applications of $LJ$ cut}
{A, \Gamma, \Gamma,\ldots,\Gamma, \Delta^{\ast} \vdash \Lambda}}
\pravilo{$LJ$ cut}}{\Bpak{\Gamma, \Gamma, \Gamma,\ldots,\Gamma, \Delta^{\ast}
\vdash \Lambda}{interchanges and contractions}{\Gamma, \Delta^{\ast}
\vdash \Lambda}}\pravilo{$LJ$ cut} \leqno{\raisebox{9ex}{($1^{\ast}$P)}} $$
Transforming (1P) into ($1^\ast$P) involves pushing cut above contraction.
The monotomic reconstruction ($1^\ast$M) of ($1^\ast$) is obtained
from (1M) by permuting interchanges with contractions, and transforming
(1M) into ($1^\ast$M) does {\em not} involve pushing cut above contraction.

The next interesting case is when we have
$$\f{\afrac{\Gamma \vdash A} \rzb \f{\Psi, \Delta \vdash \Lambda_{1}}
{A, \Delta \vdash \Lambda_{2}}\pravilo{R}}
{\Gamma, \Delta^{\ast} \vdash \Lambda_{2}}\pravilo{mix}\leqno{(2)} $$
where $A$ does not occur in $\Gamma$ and either R is introduction of
$\wedge$ on the left, in which case $A$ is of the form
$A_{1}{\wedge}A_{2}$, while $\Psi$ is either $A_{1}$ or $A_{2}$, and
$\Lambda_{1}$ is equal to $\Lambda_{2}$, or R is introduction of $\neg$
on the left, in which case $A$ is of the form ${\neg}A_{1}$, while
$\Psi$ and $\Lambda_{2}$ are empty and $\Lambda _{1}$ is $A_{1}$.
Polytomically, (2) is reconstructed as
$$ \f{\bfrac{\Gamma \vdash A}{8.3ex} \rzb \f{\afrac{\Gamma \vdash A} \rzb
\f{\Psi,\Delta \vdash \Lambda_{1}}{A, \Delta \vdash \Lambda_{2}}\pravilo{R}}
{\Bpak{\Gamma, \Delta \vdash \Lambda_{2}}{interchanges and $LJ$ cuts}
{A, \Gamma,\ldots,\Gamma, \Delta^{\ast} \vdash \Lambda_{2}}}
\pravilo{$LJ$ cut}}{\Bpak{\Gamma, \Gamma,\ldots,\Gamma, \Delta^{\ast}
\vdash \Lambda_{2}}{interchanges and contractions}{\Gamma, \Delta^{\ast}
\vdash \Lambda_{2}}}\pravilo{$LJ$ cut} \leqno{\raisebox{3ex}{(2P)}} $$
and monotomically as
$$\f{\bfrac{\Gamma \vdash A}{8.3ex} \rzb \f{\Psi, \Delta \vdash \Lambda_{1}}
{\Bpak{A, \Delta \vdash \Lambda_{2}}{interchanges and contractions}
{A, \Delta^{\ast} \vdash \Lambda_{2}}}
\pravilo{R}}{\Gamma, \Delta^{\ast} \vdash \Lambda_{2}}
\pravilo{$LJ$ cut} \leqno{\raisebox{7ex}{(2M)}} $$
In [1935] (III.3.121.22 and 3.121.222) Gentzen transforms (2) into
$$\f{\bfrac{\Gamma \vdash A}{13.5ex} \rzb
\f{\Gamma \vdash A \rzb \Psi, \Delta \vdash \Lambda_{1}}
{\f{\Bpak{\Gamma, \Psi^{\ast}, \Delta^{\ast} \vdash \Lambda_{1}}
{thinning or interchanges}{\Psi, \Gamma, \Delta^{\ast} \vdash \Lambda_{1}}}
{A, \Gamma, \Delta^{\ast} \vdash \Lambda_{2}}\pravilo{R}}
\pravilo{mix}}{\Bpak{\Gamma, \Gamma, \Delta^{\ast} \vdash \Lambda_{2}}
{interchanges and contractions}{\Gamma, \Delta^{\ast} \vdash \Lambda_{2}}}
\pravilo{mix (i.e. $LJ$ cut)} \leqno{\raisebox{7ex}{($2^{\ast}$)}} $$

When the upper mix of ($2^{\ast}$) is reconstructed polytomically, the result
of the reconstruction being called ($2^{\ast}$P), transforming (2P) into
($2^{\ast}$P) involves permuting $LJ$ cuts with $LJ$ cuts and with R.
(This permuting of cut with cut corresponds to (3.8) of the proof of
Theorem 6.1 below, and not to (2.4) and (3.7).)
It also involves pushing contraction above cut, but it does {\em not}
involve pushing cut above contraction.

When, on the other hand, the upper mix of ($2^\ast$) is reconstructed
monotomically, the result of the reconstruction being called ($2^\ast$M),
transforming (2M) into ($2^\ast$M) involves, among other things, pushing
cut above contraction.

The final interesting case is when we have
$$\f{\afrac{\Gamma \vdash A} \rzb \f
{\Delta \vdash B \rzb C, \Theta \vdash \Lambda}
{B \str C, \Delta, \Theta \vdash \Lambda}\pravilo{$\str$L}}
{\Gamma, (B \str C)^\ast, \Delta^{\ast}, \Theta^\ast \vdash \Lambda}
\pravilo{mix}\leqno{(3)} $$
where $A$ does not occur in $\Gamma$, while $(B \str C)^\ast$ stands
either for the empty sequence or for $B \str C$, according as $A$ is
$B \str C$ or not, and $A$ occurs in both $\Delta$ and $\Theta$.
Polytomically, (3) is reconstructed as 
$$ \f{\bfrac{\Gamma \vdash A}{8.3ex} \rzb \f{\bfrac{\Gamma \vdash A}{8.3ex}
\rzb \f{\Delta \vdash B \rzb C, \Theta \vdash \Lambda}
{\Bpak{B{\str}C, \Delta, \Theta \vdash \Lambda}{interchanges}
{A, A,\ldots, A, B{\str}C, \Delta^\ast, \Theta^\ast \vdash \Lambda}}
\pravilo{$\str$L}}
{\Bpak{\Gamma, A,\ldots,A, B{\str}C, \Delta^\ast, \Theta^\ast \vdash \Lambda}
{interchanges and $LJ$ cuts}
{A, \Gamma,\ldots,\Gamma, (B{\str}C)^\ast, \Delta^{\ast}, \Theta^\ast
\vdash \Lambda}}\pravilo{$LJ$ cut}}
{\Bpak{\Gamma, \Gamma,\ldots,\Gamma, (B{\str}C)^\ast, \Delta^{\ast}, \Theta^\ast
\vdash \Lambda}{interchanges and contractions}
{\Gamma, (B{\str}C)^\ast, \Delta^{\ast}, \Theta^\ast
\vdash \Lambda}}\pravilo{$LJ$ cut} \leqno{\raisebox{3ex}{(3P)}} $$
and monotomically as
$$\f{\bfrac{\Gamma \vdash A}{8.3ex} \rzb \f{\Delta \vdash B \rzb
C, \Theta \vdash \Lambda}
{\Bpak{B{\str}C, \Delta, \Theta \vdash \Lambda}{interchanges and contractions}
{A, (B{\str}C)^\ast, \Delta^{\ast}, \Theta^\ast \vdash \Lambda}}
\pravilo{$\str$L}}{\Gamma, (B{\str}C)^\ast, \Delta^{\ast}, \Theta^\ast
\vdash \Lambda}\pravilo{$LJ$ cut} \leqno{\raisebox{7ex}{(3M)}} $$
In [1935] (III.3.121.233.1) Gentzen transforms (3) into
$$\f{\bfrac{\f{\Gamma \vdash A \rzb \Delta \vdash B}
{\Gamma, \Delta^\ast \vdash B}\pravilo{mix}}{3.5ex} \rzc
\f{\Gamma \vdash A \rzb C, \Theta \vdash \Lambda}
{\Bpak{\Gamma, C^\ast, \Theta^\ast \vdash \Lambda}{thinning or interchanges}
{C, \Gamma, \Theta^\ast \vdash \Lambda}}\pravilo{mix}}
{B{\str}C, \Gamma, \Delta^\ast, \Gamma, \Theta^\ast \vdash \Lambda}
\pravilo{($\str$L)} \leqno{\raisebox{3ex}{($3^\ast$)}} $$
which if $B \str C$ is $A$, is continued by
\[\f{\Gamma \vdash A \rzb B{\str}C, \Gamma, \Delta^\ast, \Gamma, \Theta^\ast
\vdash \Lambda}
{\Bpak{\Gamma, \Gamma, \Delta^\ast, \Gamma, \Theta^\ast \vdash \Lambda}
{interchanges and contractions}
{\Gamma, \Delta^\ast, \Theta^\ast \vdash \Lambda}}
\pravilo{mix (i.e. $LJ$ cut)}\]
and if $B \str C$ is not $A$, is continued by
\[\Bpak{B{\str}C, \Gamma, \Delta^\ast, \Gamma, \Theta^\ast \vdash \Lambda}
{interchanges and contractions}
{\Gamma, B{\str}C, \Delta^\ast, \Theta^\ast \vdash \Lambda}\]
When the two top mixes of ($3^\ast$) are reconstructed polytomically,
the result of the reconstruction being called ($3^\ast$P), transforming
(3P) into ($3^\ast$P), which is analogous to the transformation of
(2P) into ($2^\ast$P), does {\em not} involve pushing cut above
contraction.  When, on the other hand, the two top mixes of ($3^\ast$)
are reconstructed monotomically, the result of the reconstruction being
called ($3^\ast$M), transforming (3M) into ($3^\ast$M) involves pushing
cut above contraction. 

So we can conclude that in the polytomic
reconstruction pushing cut above contraction is involved in the
first case, while in the monotomic reconstruction in the second and
third case. If in the first case we favour the monotomic reconstruction,
while nonuniformly, in the second and third case we favour the polytomic
reconstruction, we shall never have to push cut above contraction in order
to perform the steps of Gentzen's procedure, but we shall need this
pushing to pass from a monotomic reconstruction to the corresponding
polytomic reconstruction.

It is worth remarking that in the polytomic reconstruction, in the second and
third case we don't only lack pushing cut above contraction, but instead,
in the opposite direction, we push contraction above cut.

In the second case and in the third case when the mix formula $A$ is
$B \str C$, let us call the lowest mix in ($2^\ast$) and ($3^\ast$).
which is in fact an $LJ$ cut, the {\em critical mix} of the transformation.
The specificity of the critical mix is that it is the lowest mix in
the figure and that its right rank is 1. In the monotomic reconstruction,
the critical mix, i.e. $LJ$ cut, originates from one of the two cuts obtained
by pushing a cut above a contraction. In this pushing, which is the
relativization to $LJ$ of the transformation of ($\ast$) into ($\ast\ast$)
of Section 1, we must ensure that the critical mix originates in the
lower cut of ($\ast\ast$). Otherwise, we would need to permute also
cut with cut to ensure that the critical mix ends up as being the lowest cut.

\section{A sequent system for intuitionistic propositional logic}

\noindent Our propositional language will have the propositional constant
$\bot$ and the binary connectives $\I$, $\vee$ and $\str$. We use
$A, B, C, \ldots , A_{1}, \ldots$ as schematic letters for formulae and
$\Gamma, \Delta, \Theta, \ldots , \Gamma_{1}, \ldots$ as schematic letters
for finite, possibly empty, sequences of formulae. As usual, $\neg A$ can
be defined as $A \str \bot$. Sequents are expressions of the form
$\Gamma \vdash A$.

The sequent system \G\ has as postulates the postulates of
the sequent system \Gr\ below with all superscripts omitted.
The postulates of \G\ are named by the same names as in \Gr\
save that the superscript {\em r} is always omitted.
(It would be wasteful to write these postulates twice, once for \G\
without superscripts, and once again, just a little bit
further down, for
\Gr\, with the superscripts added.)

In Gentzen's original rules of [1935] the sequence $\Theta$ 
in the postulates of \G\ is always empty,
both in the structural rules and in the rules for connectives. Our, more
general, rules are derivable from Gentzen's rules in the presence of the
structural rule of interchange. We already replaced (Gentzen's cut)
by the present form of cut in Section 1, in order to be able
to permute cut with the structural rule of interchange. We replace
likewise the other rules of Gentzen by the present more general forms to
be able to permute contraction with other rules, and, also, for the sake
of uniformity.

As usual, we call an application of (cut) in a proof of \G\ {\em a cut}.
With this form of speech it should be kept in mind that  
our cuts are applications of the rule (cut) of \G, and not of 
(Gentzen's cut).

The {\em degree} of a cut is, as usual, the number of binary connectives in
the cut formula $A$. 
The {\em degree of a proof} in \G\ is the maximal degree among the degrees of
the cuts in this proof. A proof of degree 0 can have only
cuts whose cut formulae are atomic. A proof without cuts has degree 0.

To compute the rank of a cut we introduce an auxiliary
sequent system we call \Gr. In the sequents of \Gr\ we don't have ordinary
formulae, but indexed formulae $A^n$ where $A$ is an ordinary formula
and the {\em rank index} $n \geq 1$ is a natural number. To formulate
the postulates of \Gr\ we introduce the following conventions. If
$\Gamma^i$ is a sequence of indexed formulae, then $\Gamma^{i+1}$ is
the sequence of indexed formulae obtained by increasing by 1 every
rank index in $\Gamma^i$. (Note that
here the subscript $i$ does not stand for a single natural number: it is a
schema for any natural number in the rank
indices of the sequence $\Gamma^i$.)
We use $\Gamma^i$ and $\Gamma^j$ for sequences
of indexed formulae that may differ only in the rank indices. 
When for $\Gamma^i$ and $\Gamma^j$ we write $i \leq j$, that
means that if in $\Gamma^i$ we
find $A^n$ and at the same place in $\Gamma^j$ we find $A^m$, then
$n \leq m$. Starting from
$\Gamma^i$ and $\Gamma^j$ we obtain the sequence of indexed formulae
$\Gamma^{max(i,j)+1}$ in the following manner: if in $\Gamma^i$ we
find $A^n$ and at the same place in $\Gamma^j$ we find $A^m$, then at the
same place in $\Gamma^{max(i,j)+1}$ we put $A^{max(n,m)+1}$.

We can now give the postulates of \Gr, which are just indexed variants
of the postulates of \G:
\[
\begin{array}{llrl}
\makebox[3ex][l]{\mbox{\em axioms}} & & \Pime{structural rules} &\\
& ({\tekst{1}}^r) \rza A^1 \vdash A^1 \rzc & &
({\tekst C}^r) \rza \f{\Delta^i, A^n, B^m, \Gamma^j \vdash C^k}
{\Delta^{i+1}, B^{m+1}, A^{n+1}, \Gamma^{j+1} \vdash C^{k+1}} \\
& (\bot^r) \rza \bot^1 \vdash A^1 &
& ({\tekst W}^r) \rza \f{\Theta^i, A^n, A^m, \Gamma^j \vdash C^k}
{\Theta^{i+1}, A^{max(n,m)+1}, \Gamma^{j+1} \vdash C^{k+1}} \\
& & & ({\tekst K}^r) \rza \f{\Theta^i, \Gamma^j \vdash C^k}
{\Theta^{i+1}, A^1, \Gamma^{j+1} \vdash C^{k+1}} \\
&{\mbox{\hspace{35ex}}}
& & ({\tekst{cut}}^r) \rza \f{\Delta^i \vdash A^n \rzb \Theta^j, A^m,
\Gamma^h \vdash C^k}
{\Theta^{j+1}, \Delta^{i+1}, \Gamma^{h+1} \vdash C^{k+1}}\\[3em]
\end{array}
\]
\[
\begin{array}{llrl}
\Pime{rules for connectives} \cr \\
& ({\I\tekst{L}}^r) \rza \f{\Theta^i, A^n, \Gamma^j \vdash C^k}
{\Theta^{i+1}, {A \I B}^1, \Gamma^{j+1} \vdash C^{k+1}} &
{\mbox{\hspace{7ex}}}&
\f{\Theta^i, B^n, \Gamma^j \vdash C^k}{\Theta^{i+1}, {A \I B}^1,
\Gamma^{j+1} \vdash C^{k+1}} \\
& ({\I\tekst{R}}^r) \rza \f{\Gamma^i \vdash A^n \rzb \Gamma^j \vdash B^m}
{\Gamma^{max(i,j)+1} \vdash {A \I B}^1} & & \\
& ({\vee\tekst{L}}^r) \rza {\makebox[3ex][l]
{$\f{\Theta^i, A^n, \Gamma^j \vdash C^k
\rzb \Theta^l, B^m, \Gamma^h \vdash C^q}
{\Theta^{max(i,l)+1}, {A \vee B}^1, \Gamma^{max(j,h)+1}
\vdash C^{max(k,q)+1}}$}} & &
\\
& ({\vee\tekst{R}}^r) \rza \f{\Gamma^i \vdash A^n}{\Gamma^{i+1}
\vdash {A \vee B}^1} & &
\f{\Gamma^i \vdash B^n}{\Gamma^{i+1} \vdash {A \vee B}^1} \\
& ({\str\tekst{L}}^r) \rza \f{\Delta^i \vdash A^n \rzb
\Theta^j, B^m, \Gamma^h \vdash C^k}
{\Theta^{j+1}, \Delta^{i+1}, {A \str B}^1, \Gamma^{h+1} \vdash C^{k+1}}
& & ({\str\tekst{R}}^r) \rza \f{A^n, \Gamma^i \vdash B^m}
{\Gamma^{i+1} \vdash {A \str B}^1}
\end{array} \]

Take a cut in \G, and do again in \Gr\ the proofs of the two premises
$\Delta \vdash A$ and $\Theta, A, \Gamma \vdash C$ of this cut exactly
as they are done in \G, save that in \Gr\ rank indices are taken into account.
Let these two proofs in \Gr\ prove $\Delta^i \vdash A^n$ and $\Theta^j,
A^m, \Gamma^h \vdash C^k$. Then the {\em rank} of our cut is
$n+m$. The {\em left rank} of this cut is $n$, and the {\em right rank}
is $m$.

A cut in \G\ is {\em topmost} iff there are no cuts above it. Gentzen 
computed rank only for topmost mixes, i.e. those above which
there are no mixes, and his notion of rank coincides with our
notion of rank for topmost cuts. We need, however, for our
cut-elimination procedure (see Section 6 below) 
the more general notion of rank
we have just introduced, which applies to
any cut, and not only topmost cuts.

\section{Cut elimination in implicationless \G}

\noindent For every proof in \G\ in which the connective of implication
$\str$ does not occur, there is a simple procedure of cut elimination, which
eliminates cut directly, not via mix, and involves pushing cut above
contraction. To describe this procedure we introduce the following auxiliary
implicationless sequent system called \Gz. (The index $z$ stands for
``Zucker'', from whose indexing of sequents in [1974], the indices of \Gz\
are derived; a measure analogous to these indices may be found in
[B. 1997].)

On the left-hand sides of the sequents of \Gz\ we don't have ordinary formulae,
but indexed formulae $A^\alpha$ where $A$ is an ordinary implicationless
formula and the {\em contraction index} $\alpha \geq 1$ is a natural number.
To formulate the postulates of \Gz\ we use conventions analogous to those
we used for \Gr\ in the preceding section.

The postulates of \Gz\ are the following indexed variants
of the postulates of \G\ minus the rules for implication:
\[
\begin{array}{llrl}
\makebox[3ex][l]{\mbox{\em axioms}} & & \Pime{structural rules} &\\
& ({\tekst{1}}^z) \rza A^1 \vdash A \rzc & &
({\tekst C}^z) \rza \f{\Delta^i, A^\alpha, B^\beta, \Gamma^j \vdash C}
{\Delta^{i}, B^{\beta}, A^{\alpha}, \Gamma^{j} \vdash C} \\
& (\bot^z) \rza \bot^1 \vdash A &
& ({\tekst W}^z) \rza \f{\Theta^i, A^\alpha, A^\beta, \Gamma^j \vdash C}
{\Theta^{i}, A^{\alpha + \beta}, \Gamma^{j} \vdash C} \\
& & & ({\tekst K}^z) \rza \f{\Theta^i, \Gamma^j \vdash C}
{\Theta^{i}, A^1, \Gamma^{j} \vdash C}\\
& {\mbox{\hspace{35ex}}} & & ({\tekst{cut}}^z)
\rza \f{\Delta ^i \vdash A \rzb \Theta^j,
A^\alpha, \Gamma^h \vdash C}
{\Theta^j, \Delta^{i\alpha}, \Gamma^h \vdash C}\\[3em]
\end{array}
\]
\[
\begin{array}{llrl}
\Pime{rules for connectives} \cr \\
& ({\I\tekst{L}}^z) \rza \f{\Theta^i, A^\alpha, \Gamma^j \vdash C}
{\Theta^{i}, {A \I B}^\alpha, \Gamma^{j} \vdash C} &
{\mbox{\hspace{7ex}}} &
\f{\Theta^i, B^\alpha, \Gamma^j \vdash C}{\Theta^{i}, {A \I B}^\alpha,
\Gamma^{j} \vdash C} \\
& ({\I\tekst{R}}^z) \rza \f{\Gamma^i \vdash A \rzb \Gamma^j \vdash B}
{\Gamma^{max(i,j)} \vdash {A \I B}} & & \\
& ({\vee\tekst{L}}^z) \rza {\makebox[35ex][l]
{$\f{\Theta^i, A^\alpha, \Gamma^j \vdash C
\rzb \Theta^l, B^\beta, \Gamma^h \vdash C}
{\Theta^{max(i,l)}, {A \vee B}^{max(\alpha, \beta)}, \Gamma^{max(j,h)}
\vdash C}$}} & &
\\
& ({\vee\tekst{R}}^z) \rza \f{\Gamma^i \vdash A}{\Gamma^{i}
\vdash {A \vee B}} & &
\f{\Gamma^i \vdash B}{\Gamma^{i} \vdash {A \vee B}} 
\end{array} \]

As we did for applications of (cut) in \G, we 
now call {\em cuts} applications of (${\tekst{cut}}^z$)
in \Gz.

We shall prove the following theorem by eliminating cut directly
and by pushing cut above contractions.
\\[0.3cm]
{\bf Theorem 4.1}\hspace{1em}
{\em  
Every proof of $\Pi^i \vdash C$ in \Gz\ can be reduced to a cut-free proof of
$\Pi^j \vdash C$ where $j \leq i$.}
\\[0.2cm]
\dkz We proceed by an induction on triples $\langle d, z, r \rangle$,
lexicographically ordered, where $d$ is the degree of a cut, $z$ is the
contraction index of the cut formula in the right premise of
(${\tekst{cut}}^z$) and {\em r} is 
the rank of the cut (rank is defined for \Gz\ as it
is defined for \G, via \Gr). We show 
that every proof of $\Pi^i \vdash C$ with
a single cut, which is the last rule of the proof, can be reduced to a
cut-free proof of $\Pi^j \vdash C$ where $j \leq i$.
\\[0.1cm]
(1)\rza Suppose the rank of our cut is 2. Then our cut is covered by at least
one of the following cases
$$ \f{\afrac{A^1 \vdash A} \rzb \fp{\pi}{\Theta^j, A^\alpha, \Gamma^h \vdash C}}
{\Theta^j, A^\alpha, \Gamma^h \vdash C}\pravilo{$\tekst{cut}^z$}
\leqno{(1.1)}$$

Then we replace this proof by the cut-free proof $\pi$ of the right premise of
$\tekst{cut}^z$.
$$ \f{\afrac{\bot^1 \vdash A} \rzb \fp{\pi}{\Theta^j, A^\alpha, \Gamma^h
\vdash C}}
{\Theta^j, \bot^\alpha, \Gamma^h \vdash C}\pravilo{$\tekst{cut}^z$}
\leqno{(1.2)}$$
Then we replace this proof by
\[ \Bpak{\bot^1 \vdash C}{applications \\ of ($\tekst{K}^z$)}
{\Theta^1, \bot^1, \Gamma^1 \vdash C}\]
$$ \f{\fp{\pi}{\Gamma^i \vdash C} \rzb \afrac{C^1 \vdash C}}
{\Gamma^i \vdash C}\pravilo{$\tekst{cut}^z$}
\leqno{(1.3)}$$
Then we replace this proof by the cut-free proof $\pi$ of the left premise
of $\tekst{cut}^z$.
$$ \f{\fp{\pi_{1}}{\Gamma^i \vdash A} \rzb \f{\fp{\pi_{2}}
{\Theta^h, \Delta^j \vdash C}}{\Theta^h, A^1, \Delta^j \vdash C}
\pravilo{$\tekst{K}^z$}}
{\Theta^h, \Gamma^i, \Delta^j \vdash C}\pravilo{$\tekst{cut}^z$}
\leqno{(1.4)}$$
Then we replace this proof by the following proof
\[ \Apak{\pi_{2}}{\Theta^h, \Delta^j \vdash C}{applications \\
of ($\tekst{K}^z$)}
{\Theta^h, \Gamma^1, \Delta^j \vdash C}\]
$$ \f{\f{\fp{\pi_1}{\Gamma^i \vdash A} \rzb \fp{\pi_2}{\Gamma^j \vdash B}}
{\Gamma^{max(i,j)} \vdash A \I B}\pravilo{$\I\tekst{R}^z$} \rzc
\f{\fp{\pi}{\Theta^h, A^\alpha, \Delta^l \vdash C}}
{\Theta^h, {A \I B}^\alpha, \Delta^l \vdash C}\pravilo{$\I\tekst{L}^z$}}
{\Theta^h, \Gamma^{max(i,j)\alpha}, \Delta^l \vdash C}
\pravilo{$\tekst{cut}^z$} \leqno{(1.5)}$$
Then we replace this proof by
\[\f{\fp{\pi_1}{\Gamma^i \vdash A} \rzb
\fp{\pi}{\Theta^h, A^{\alpha}, \Delta^l \vdash C}}
{\Theta^h, \Gamma^{i\alpha}, \Delta^l \vdash C}\pravilo{$\tekst{cut}^z$}\]
We proceed analogously when $A^\alpha$ is replaced by $B^\alpha$.
$$ \f
{\f
{\fp{\pi}{\Gamma^i \vdash A}}{\Gamma^{i} \vdash A \vee B}
\pravilo{$\vee\tekst{R}^z$} \rzc
\f{\fp{\pi_1}{\Theta^h, A^\alpha, \Delta^l \vdash C} \rzb
\fp{\pi_2}{\Theta^u, B^\beta, \Delta^v \vdash C}}
{\Theta^{max(h,u)}, {A \vee B}^{max(\alpha,\beta)}, \Delta^{max(l,v)}
\vdash C}\pravilo{$\vee\tekst{L}^z$}}
{\Theta^{max(h,u)}, \Gamma^{i\;  max(\alpha,\beta)}, \Delta^{max(l,v)} \vdash C}
\pravilo{$\tekst{cut}^z$} \leqno{(1.6)} $$
Then we replace this proof by
\[\f{\fp{\pi}{\Gamma^i \vdash A} \rzb
\fp{\pi_1}{\Theta^h, A^{\alpha}, \Delta^l \vdash C}}
{\Theta^h, \Gamma^{i\alpha}, \Delta^l \vdash C}\pravilo{$\tekst{cut}^z$}\]
We proceed analogously when $\pi$ ends with $\Gamma^i \vdash B$.
\\[0.2cm]
(2) \rza Suppose the left rank of our cut is greater than 1. Then we have
the following cases.
$$\f{\f{\fp{\pi}{\Delta^i \vdash D}}{\Phi^k \vdash D}\pravilo{R}
\rzc \fp{\pi_2}{\Theta^j, D^\gamma, \Gamma^h \vdash C}}
{\Theta^j, \Phi^{k\gamma}, \Gamma^h \vdash C}\pravilo{$\tekst{cut}^z$}
\leqno{(2.1)} $$
where R is $\tekst{C}^z$, \naz{W}, \naz{K}\ or $\I$\naz{L}. Then we
replace this proof by
\[\f{\fp{\pi}{\Delta^i \vdash D} \rzc
\afrac{\Theta^j, D^\gamma, \Gamma^h \vdash C}}
{\f{\Theta^j, \Delta^{i\gamma}, \Gamma^h \vdash C}
{\Theta^j, \Phi^l, \Gamma^h \vdash C}\pravilo{R}}\pravilo{\naz{cut}}\]
When R is \naz{C}, \naz{W}\ or $\I$\naz{L}, then $l=k\gamma$, and when R
is \naz{K}, then one index $\gamma$ of $\Phi^{k\gamma}$ is replaced by
1 in $\Phi^l$.
$$ \f{\f{\fp{\pi_1}{\Theta^i, A^\alpha, \Gamma^j \vdash D} \rzb
\fp{\pi_2}{\Theta^l, B^\beta, \Gamma^h \vdash D}}
{\Theta^{max(i,l)}, {A \vee B}^{max(\alpha,\beta)}, \Gamma^{max(j,h)}
\vdash D}\pravilo{$\vee\tekst{L}^z$} \rzc
\fp{\pi}{\Delta^u, D^\gamma, \Xi^v \vdash C}}
{\Delta^u, \Theta^{max(i,l)\gamma}, {A \vee B}^{max(\alpha,\beta)\gamma},
\Gamma^{max(j,h)\gamma}, \Xi^v \vdash C}
\pravilo{$\tekst{cut}^z$} \leqno{(2.2)}$$
Then we replace this proof by
\[ \f{
\f{\fp{\pi_1}{\Theta^i, A^\alpha, \Gamma^j \vdash D} \rzb
\fp{\pi}{\Delta^u, D^\gamma, \Xi^v \vdash C}}
{\Delta^u, \Theta^{i\gamma}, A^{\alpha\gamma}, \Gamma^{j\gamma}, \Xi^v
\vdash C}\pravilo{\naz{cut}} \rzc
\f{\fp{\pi_2}{\Theta^l, B^\beta, \Gamma^h \vdash D} \rzb
\fp{\pi}{\Delta^u, D^\gamma, \Xi^v \vdash C}}
{\Delta^u, \Theta^{l\gamma}, B^{\beta\gamma}, \Gamma^{h\gamma}, \Xi^v
\vdash C}\pravilo{\naz{cut}} 
}
{\Delta^{max(u,u)}, \Theta^{max(i\gamma,l\gamma)},
{A \vee B}^{max(\alpha\gamma,\beta\gamma)},
\Gamma^{max(j\gamma,h\gamma)}, \Xi^{max(v,v)} \vdash C}
\pravilo{$\vee$\naz{L}} \]
(3) \rza Suppose the right rank of our cut is greater than 1. Then we
have the following cases.
$$ \f{\fp{\pi_{1}}{\Delta^i \vdash D} \rzb \f{\fp{\pi_{2}}
{\Theta^j, D^\gamma, \Gamma^h \vdash E}}{\Phi^l, D^\gamma, \Xi^k \vdash C}
\pravilo{R}}
{\Phi^l, \Delta^{i\gamma}, \Xi^k \vdash C}\pravilo{$\tekst{cut}^z$}
\leqno{(3.1)}$$
where R is \naz{C}, \naz{W}, \naz{K}, $\I$\naz{L} or $\vee$\naz{R}.
Then we replace this proof by
\[\f{\fp{\pi_1}{\Delta^i \vdash D} \rzc
\fp{\pi_2}{\Theta^j, D^\gamma, \Gamma^h \vdash E}}
{\f{\Theta^j, \Delta^{i\gamma}, \Gamma^h \vdash E}
{\Phi^l, \Delta^{i\gamma}, \Xi^k \vdash C}\pravilo{R}}\pravilo{\naz{cut}}\]
save when R is \naz{C}, and when in the transformed proof R can be a number of
applications of (\naz{C}).
$$ \f
{{\fp{\pi}{\Delta^i \vdash D}} \rzc
\f{\fp{\pi_1}{\Theta^j, D^\alpha, \Gamma^h \vdash C_1} \rzb
\fp{\pi_2}{\Theta^k, D^\beta, \Gamma^l \vdash C_2}}
{\Theta^{max(j,k)}, {D}^{max(\alpha,\beta)}, \Gamma^{max(h,l)}
\vdash C_1 \I C_2}\pravilo{$\I\tekst{R}^z$}}
{\Theta^{max(j,k)}, \Delta^{i\;  max(\alpha,\beta)}, \Gamma^{max(h,l)} \vdash
C_1 \I C_2}
\pravilo{$\tekst{cut}^z$} \leqno{(3.2)} $$
Then we replace this proof by
\[ \f{
\f{\fp{\pi}{\Delta^i \vdash D} \rzb
\fp{\pi_1}{\Theta^j, D^\alpha, \Gamma^h \vdash C_1}}
{\Theta^j, \Delta^{i\alpha}, \Gamma^{h} \vdash C_1}\pravilo{\naz{cut}} \rzc
\f{\fp{\pi}{\Delta^i \vdash D} \rzb
\fp{\pi_2}{\Theta^k, D^\beta, \Gamma^l \vdash C_2}}
{\Theta^{k}, \Delta^{i\beta}, \Gamma^{l} \vdash C_2}\pravilo{\naz{cut}} 
}
{\Theta^{max(j,k)}, \Delta^{max(i\alpha,i\beta)},
\Gamma^{max(h,l)} \vdash C_1 \I C_2}
\pravilo{$\I$\naz{R}} \]
$$ \f
{{\fp{\pi}{\Delta^i \vdash D}} \rzb
\f{\fp{\pi_1}{\Theta^{j_1}_1, A^\gamma, \Theta^{j_2}_2, D^\alpha,
\Gamma^h \vdash C} \rzb
\fp{\pi_2}{\Theta^{k_1}_1, B^\delta, \Theta^{k_2}_2, D^\beta, \Gamma^l \vdash C}}
{\Theta^{max(j_1,k_1)}_1, {A \vee B}^{max(\gamma,\delta)},
\Theta^{max(j_2,k_2)}_2, D^{max(\alpha,\beta)}, \Gamma^{max(h,l)}
\vdash C}\pravilo{$\vee\tekst{L}^z$}}
{\Theta^{max(j_1,k_1)}_1, {A \vee B}^{max(\gamma,\delta)},
\Theta^{max(j_2,k_2)}_2, \Delta^{i\;  max(\alpha,\beta)}, \Gamma^{max(h,l)}
\vdash C}\pravilo{$\tekst{cut}^z$} \leqno{(3.3)} $$
Then we replace this proof by
\[ \f{
\f{\fp{\pi}{\Delta^i \vdash D} \rzb
\fp{\pi_1}{\Theta^{j_1}_1, A^\gamma, \Theta^{j_2}_2, D^\alpha,
\Gamma^h \vdash C}}
{\Theta^{j_1}_1, A^\gamma, \Theta^{j_2}_2, \Delta^{i\alpha}, \Gamma^{h}
\vdash C}\pravilo{\naz{cut}} \rzc
\f{\fp{\pi}{\Delta^i \vdash D} \rzb
\fp{\pi_2}{\Theta^{k_1}_1, B^\delta, \Theta^{k_2}_2, D^\beta, \Gamma^l
\vdash C}}
{\Theta^{k_1}_1, B^\delta, \Theta^{k_2}_2, \Delta^{i\beta}, \Gamma^{l}
\vdash C}\pravilo{\naz{cut}}
}
{\Theta^{max(j_1,k_1)}_1, {A \vee B}^{max(\gamma,\delta)},
\Theta^{max(j_2,k_2)}_2, \Delta^{max(i\alpha,i\beta)},
\Gamma^{max(h,l)} \vdash C}
\pravilo{$\vee$\naz{L}} \]
We proceed analogously when $\pi_1$ ends with
$\Theta^j, D^\alpha, \Gamma^{h_1}_1, A^\gamma, \Gamma^{h_2}_2 \vdash C$
and $\pi_2$ ends with
$\Theta^k, D^\beta, \Gamma^{l_1}_1, B^\delta, \Gamma^{l_2}_2 \vdash C$.
$$ \f{\fp{\pi_{1}}{\Delta^i \vdash D} \rzb \f{\fp{\pi_{2}}
{\Theta^j, D^\alpha, D^\beta, \Gamma^h \vdash C}}{\Theta^j, D^{\alpha+\beta},
\Gamma^h \vdash C}
\pravilo{\naz{W}}}
{\Theta^j, \Delta^{i(\alpha+\beta)}, \Gamma^h \vdash C}\pravilo{$\tekst{cut}^z$}
\leqno{(3.4)}$$
Then we replace this proof by
\[\f{\fp{\pi_{1}}{\Delta^i \vdash D} \rzb \f{\fp{\pi_1}{\Delta_i \vdash D}
\rzb \fp{\pi_{2}}{\Theta^j, D^\alpha, D^\beta, \Gamma^h \vdash C}}
{\Theta^j, \Delta^{i\alpha}, D^{\beta}, \Gamma^h \vdash C}
\pravilo{\naz{cut}}}
{\Bpak{\Theta^j, \Delta^{i\alpha}, \Delta^{i\beta}, \Gamma^h \vdash C}
{applications \\ of (\naz{C}) and (\naz{W})}{\Theta^j, \Delta^{i\alpha+i\beta},
\Gamma^h \vdash C}}\pravilo{$\tekst{cut}^z$} \]
In the transformed proof, the cut formula of the upper \naz{cut}\
has the same degree as in the original cut, but it has a lower contraction
index in the right premise (even the rank has decreased, but this is not now
essential). Hence by the induction hypotheses we have a cut-free proof of
$\Theta^k, \Delta^l, D^\gamma, \Gamma^n \vdash C$ with $k \leq j$,
$l \leq i\alpha$, $\gamma \leq \beta$ and $n \leq h$. So we obtain
\[\f{\fp{\pi_{1}}{\Delta^i \vdash D} \rzc
\fp{\pi}{\Theta^k, \Delta^l, D^\gamma, \Gamma^n \vdash C}}
{\Bpak{\Theta^k, \Delta^{l}, \Delta^{i\gamma}, \Gamma^n \vdash C}
{applications \\ of (\naz{C}) and (\naz{W})}{\Theta^k, \Delta^{l+i\gamma},
\Gamma^n \vdash C}}\pravilo{$\tekst{cut}^z$} \]
where the cut formula of \naz{cut}\ is again of the same degree as
in the original cut, but has a lower contraction index in the right premise
(its rank has perhaps increased).
\qed

The contraction indices of \Gz\ are not the only possible indices we could
have chosen. For example, we could replace $\alpha + \beta$ by
$max(\alpha,\beta) + 1$ in (\naz{W}). Whereas the original contraction index
measures the number of contractions in the clusters, this new 
index would measure the height of clusters. (For the notion of
cluster, in German {\em Bund}, 
see [Gentzen 1938], Section 3.41; see also [D. \& P. 1999] and
references therein.) A rationale for the maximum function in the indices of
$\Gamma$ and $\Theta$ in ($\I$\naz{R}) and ($\vee$\naz{L}) may be found in
the proofs of Lemma 5.3 and Theorem 5.5 below.

\section{W-normal form}

\noindent To formulate our new cut-elimination procedure for \G\ we need to
introduce the following notion of normal form. (Note that W is sometimes used 
as a label for thinning, also called "weakening", while our
use of this label for contraction is suggested by combinatory logic.
So our terminology shouldn't be confused with the terminology of some
other authors, which may use the same terms to designate other
things; cf., for example, [Mints 1996].)

A proof in \G\ is called W-{\em normal} iff every application of (W)
in this proof is either the last rule of the proof, or it has only
applications of (W) below it, or it is the upper rule in the following
contexts:
\[\f{\Theta, A, A, A, \Gamma \vdash C}{\f{\Theta, A, A, \Gamma \vdash C}
{\Theta, A, \Gamma \vdash C}\pravilo{W}}\pravilo{W}\]
\[\f{A, A, \Gamma \vdash B}{\f{A, \Gamma \vdash B}{\Gamma \vdash A \str B}
\pravilo{$\str$R}}\pravilo{W}\]
We also need the following terminology.

We say that an application of (W) in a proof
\[ \fp{\pi}{\Gamma \vdash C}\]
is {\em tied} to an occurrence $G$ of a formula in $\Gamma$ iff the
principal formula (i.e. contracted formula) of this application of
(W) belongs to the cluster of $G$ in $\pi$. (For the notion of cluster see
[Gentzen 1938], Section 3.41.)

A contraction in a proof $\pi$ is {\em engaged} iff it is tied to the cut
formula of the right premise of some cut in $\pi$. If the corresponding cut is
immediately below the engaged contraction, then we call such a contraction
{\em directly engaged}. A contraction in $\pi$ that is not engaged  is
called {\em neutral}.

We shall prove now a series of lemmata leading to the proof of the theorem
that every proof can be reduced to a W-normal proof of the same degree
of the same sequent. This theorem covers the first phase of our
cut-elimination procedure.
\\[0.3cm]
\lema{5.1}
{\em Every segment of a proof $\pi$ of the form
\[ \Dpak{\Phi \vdash C}{$e+n$ applications of (W) followed by applications of
(C)}{\Psi \vdash C}\]
can be transformed into a segment of the form
\[ \Dpak{\Phi \vdash C}{applications of (C) followed by $e+n$ applications
of (W)}{\Psi \vdash C}\]
where $e$ is the number of engaged contractions of $\pi$ and $n$ the number of
neutral contractions of $\pi$ that occur in the figures above. The
degree of the transformed proof is the same as the degree of $\pi$.}
\\[0.2cm]
\dkz
By induction on the lexicographically ordered couples $\langle e+n, i \rangle$,
where $i$ is the number of applications of (C) in the initial segment $\cal S$.
\\[0.3cm]
(a)\hspace{1.5ex}The segment $\cal S$ is of the form
\[\f{\Dpak{\Phi \vdash C}{$e+n-1$ applications of (W)}
{\Gamma, A, A, D, \Delta \vdash C}}{\f{\Gamma, A, D, \Delta \vdash C}
{\Dpak{\Gamma, D, A, \Delta \vdash C}{$i-1$ applications of (C)}
{\Psi \vdash C}}\pravilo{C}}\pravilo{W}\]
By transforming the segment beginning
with $\Gamma, A, A, \Delta \vdash C$ and ending with
$\Gamma, D, A, \Delta \vdash C$, the whole segment $\cal S$ is
transformed into
\[\f{\Dpak{\Phi \vdash C}{$e+n-1$ applications of (W)}
{\Gamma, A, A, D, \Delta \vdash C}}{\f{
\f{\Gamma, A, D, A, \Delta \vdash C}{\Gamma, D, A, A, \Delta \vdash C}
\pravilo{C}}
{\Dpak{\Gamma, D, A, \Delta \vdash C}{$i-1$ applications of (C)}
{\Psi \vdash C}}\pravilo{W}}\pravilo{C}\]
This transformation preserves the engagement or neutrality of the
lowest contraction.

Since $\langle e+n-1, 
2 \rangle < \langle e+n, i
\rangle$, by the induction 
hypothesis the segment beginning with $\Phi \vdash C$ and ending with
$\Gamma, D, A, A, \Delta \vdash C$ can be
transformed so that our whole segment becomes
\[\f{\Dpak{\Dpak{\Phi \vdash C}{applications of (C)}{\Phi' \vdash C}}
{$e+n-1$ applications of (W)}
{\Gamma, D, A, A, \Delta \vdash C}}{\Dpak{\Gamma, D, A, \Delta \vdash C}
{$i-1$ applications of (C)}
{\Psi \vdash C}}\pravilo{W}\]

Since $\langle e+n, i-1 \rangle < \langle e+n, i \rangle $, 
by the induction hypothesis the segment
beginning with $\Phi' \vdash C$ and ending with $\Psi \vdash C$ 
can be transformed so that the whole segment 
is brought into the desired form.
\\[0.3cm]
(b)\hspace{1.5ex}The segment $\cal S$ is of the form
\[\f{\Dpak{\Phi \vdash C}{$e+n-1$ applications of (W)}
{\Gamma, A, A, \Delta \vdash C}}{\f{\Gamma, A, \Delta \vdash C}
{\Dpak{\Gamma', A, \Delta' \vdash C}{$i-1$ applications of (C)}
{\Psi \vdash C}}\pravilo{C}}\pravilo{W}\]
By transforming the 
segment beginning with $\Gamma, A, A, \Delta \vdash C$ and ending
with $\Gamma', A, \Delta' \vdash C$, the whole segment $\cal S$ 
is transformed into
\[\f{\Dpak{\Phi \vdash C}{$e+n-1$ applications of (W)}
{\Gamma, A, A, \Delta \vdash C}}{\f{
\Gamma', A, A, \Delta' \vdash C}
{\Dpak{\Gamma', A, \Delta' \vdash C}{$i-1$ applications of (C)}
{\Psi \vdash C}}\pravilo{W}}\pravilo{C}\]
The remaining steps are analogous to the steps in (a).
\qed

A W-normal proof is called {\em tailless} iff its last rule is not an
application of (W).
Let $\pi_1$ and $\pi_2$ be tailless. We define inductively as follows
the class of proofs ${\cal C}(\pi_1, \pi_2)$:
\\[0.3cm]
(i)\hspace{2ex}The proof $\pi_2$ belongs to ${\cal C}(\pi_1, \pi_2)$.
\\[0.2cm]
(ii)\hspace{1.5ex}If $\pi$ belongs to ${\cal C}(\pi_1, \pi_2)$, 
then the proof
\[\Cpak{\f{\pi_1 \rzb \pi}{\Phi \vdash B}\pravilo{cut}}{applications of (C)}
{\Psi \vdash B}\]
belongs to ${\cal C}(\pi_1, \pi_2)$,
provided that there is no occurrence of a formula in a subproof
$\pi_1$ of $\pi$ that belongs to the cluster of the cut formula
in the right premise of the cut noted in the figure.
\\[0.2cm]
(iii)\hspace{1ex}If $\pi$ belongs to ${\cal C}(\pi_1, \pi_2)$, 
then $\pi$ followed by an
application of (W) belongs to ${\cal C}(\pi_1, \pi_2)$.
\\[0.3cm]  

The application of (W) in (iii) in the definition of ${\cal C}(\pi_1,\pi_2)$
is called {\em mobile}.
The {\em height} of a mobile application of (W) is the number of
applications of (W) and (cut) below it in the proof (we don't count
applications of (C)).

It is easy to verify that an application of (W) in a tailless subproof
of a proof cannot be engaged in this proof. This fact will be useful in the
proof of the following lemma.
\\[0.3cm]
\lema{5.2}
{\em For every pair of tailless proofs $\pi_1$ and $\pi_2$, every proof
from ${\cal C}(\pi_1,\pi_2)$ can be transformed into a W-normal proof
of the same degree.}
\\[0.2cm]
\dkz
By induction on the lexicographically ordered pairs $\langle \kappa,
\lambda \rangle$, where $\kappa$ is the number of engaged
applications of (W) in the proof and $\lambda$ is the sum of the heights of all
mobile applications of (W) in the proof.

By the definition of \cls, if there is no mobile application of (W)
followed immediately by a cut in a proof from \cls, then
this proof is W-normal.
\\[0.3cm]
(a)\hspace{1.5ex}Suppose our proof is of the form
\[\f{\fp{\pi_1}{\Gamma \vdash C} \rzb \f{\fp{\pi}
{\Theta, C, C, \Delta \vdash B}}{\Theta, C, \Delta \vdash B}
\pravilo{W directly engaged}}
{\Dpak{\Theta, \Gamma, \Delta \vdash B}{applications of (C) followed by \\
applications of (W)}{\fp{\Xi \vdash B}{\vdots}}}\pravilo{cut}\]
for $\pi$ in \cls.

By pushing the directly engaged application of (W) below 
cut, this proof is transformed into 
\[\f{\fp{\pi_1}{\Gamma \vdash C} \rzb \f{\fp{\pi_1}{\Gamma \vdash C}
\rzb \fp{\pi}
{\Theta, C, C, \Delta \vdash B}}{\Theta, \Gamma, C, \Delta \vdash B}
\pravilo{cut}}
{\Dpak{\Dpak{\Theta, \Gamma, \Gamma, \Delta \vdash B}
{applications of (C) followed by \\ neutral applications of (W)}
{\Theta, \Gamma, \Delta \vdash B}}{remaining applications of (C) \\
followed by applications of (W)}{\fp{\Xi \vdash B}{\vdots}}}\pravilo{cut}\]
The neutral applications of (W) mentioned above, which contract formulae
from $\Gamma$, are neutral by the proviso in (ii) of the definition of \cls.

By Lemma 5.1, the segment beginning with
$\Theta, \Gamma, \Gamma, \Delta \vdash B$ and ending with $\Xi \vdash B$
can be transformed so that our whole proof becomes
\[\f{\fp{\pi_1}{\Gamma \vdash C} \rzb \f{\fp{\pi_1}{\Gamma \vdash C}
\rzb \fp{\pi}
{\Theta, C, C, \Delta \vdash B}}{\Theta, \Gamma, C, \Delta \vdash B}
\pravilo{cut}}
{\Dpak{\Theta, \Gamma, \Gamma, \Delta \vdash B}
{applications of (C) followed by \\applications of (W)}
{\fp{\Xi \vdash B}{\vdots}}}\pravilo{cut}\]
which belongs to \cls\ and has one engaged application of (W) less than the
original proof. (Here we use the fact that no application of (W) in
$\pi_1$ can be engaged in the proof above.) So the measure of the transformed
proof is $\langle \kappa-1, \lambda \rangle < \langle \kappa, \lambda \rangle$.
By the induction hypothesis this proof can be transformed into a W-normal
proof.
\\[0.3cm]
(b)\hspace{1.5ex}Suppose our proof is of the form
\[\f{\fp{\pi_1}{\Gamma \vdash C} \rzb \f{\fp{\pi}
{\Theta, C, \Delta \vdash B}}{\Theta', C, \Delta' \vdash B}
\pravilo{W neutral or not directly engaged}}
{\Dpak{\Theta', \Gamma, \Delta' \vdash B}
{applications of (C) followed by \\ applications of (W)}
{\fp{\Xi \vdash B}{\vdots}}}\pravilo{cut}\]
By pushing the distinguished application of (W), which immediately
follows $\pi$, below cut, this proof is transformed into
\[\f{\fp{\pi_1}{\Gamma \vdash C} \rzb \fp{\pi}{\Theta, C, \Delta \vdash B}}
{\Dpak{\f{\Theta, \Gamma, \Delta \vdash B}
{\Theta', \Gamma, \Delta' \vdash B}\pravilo{W}}
{applications of (C) followed by \\ applications of (W)}
{\fp{\Xi \vdash B}{\vdots}}}\pravilo{cut}\]
where the distinguished applications of (W) in the original figure and in
the transformed figure are either both neutral or both engaged. By Lemma
5.1, the segment beginning with $\Theta, \Gamma, \Delta \vdash B$
and ending with $\Xi \vdash B$ can be
transformed so that our whole proof becomes
\[\f{\fp{\pi_1}{\Gamma \vdash C} \rzb \fp{\pi}{\Theta, C, \Delta \vdash B}}
{\Dpak{\Theta, \Gamma, \Delta \vdash B}
{applications of (C) followed by \\ applications of (W)}
{\fp{\Xi \vdash B}{\vdots}}}\pravilo{cut}\]
This proof belongs to \cls\ and has the same number of engaged applications
of (W), but its $\lambda$ has decreased by 1. Since 
$\langle \kappa, \lambda -1 \rangle <
\langle \kappa, \lambda \rangle$, by the induction hypothesis 
our proof can be transformed into
a W-normal proof.
\qed

The following lemma is covered by Lemma 12 in [Kleene 1952]. However,
Kleene's sequent system is not quite the same: interchange is only
implicit in it, and his proof doesn't cover all details we need
to cover. (In his proof on p. 24, in the third illustration, Kleene
assumes that the $n_1$ contractions above $A, A, \Gamma \str \Theta$
are all tied to the first $A$, whereas we cannot assume that. We
could assume it only after introducing a new reduction step that
transforms sequences of contractions tied to the same occurrence 
of a formula.)
\\[0.3cm]
\lema{5.3}
{\em A proof of the form
\[\f{\Apak{\pi_1}{\Phi \vdash C}{applications \\ of (W)}
{\Theta, A, \Gamma \vdash C} \rzc \Apak{\pi_2}{\Psi \vdash C}
{applications \\ of (W)}{\Theta, B, \Gamma \vdash C}}
{\Theta, A \vee B, \Gamma \vdash C}\pravilo{$\vee$L}\]
where $\pi_1$ and $\pi_2$ are tailless, can be transformed into a
W-normal proof, of the same degree, of the form
\[\Dpak{\fp{\pi}{\Xi \vdash C}}{applications of (W)}
{\Theta, A \vee B, \Gamma \vdash C}\]
where $\pi$ is tailless, and for every occurrence $G$ of a formula in
$\Theta$, if above the left premise of} $\vee$L {\em in the former figure
there are $k_1$ applications of} (W) {\em tied to $G$, and if above the
right premise of} $\vee$L {\em in the former figure there are $k_2$
applications of} (W) {\em tied to this 
same $G$, then in the latter figure there
are $max(k_{1}, k_2)$ applications of} (W) {\em tied to $G$. The same holds
for occurrences of formulae in $\Gamma$.}
\\[0.2cm]
\dkz
Let $n$ be the number of applications of (W) tied to $A$ in the left premise
of $\vee$L and $m$ be the number of applications of (W) tied to $B$ in the
right premise of $\vee$L in the figure of the initial proof. We prove
the lemma by induction on $n+m$. Our proof is first transformed into
\[\f{
\Bpak{\Apak{\pi_1}{\Phi \vdash C}{applications \\ of (K)}
{\Theta', A,\ldots,A, \Gamma' \vdash C}}{$n$ applications \\ of (W)}
{\Theta', A, \Gamma' \vdash C}
\rzd
\Bpak{\Apak{\pi_2}{\Psi \vdash C}{applications \\ of (K)}
{\Theta', B,\ldots,B, \Gamma' \vdash C}}{$m$ applications \\ of (W)}
{\Theta', B, \Gamma' \vdash C}
}
{\Cpak{\Theta', A \vee B, \Gamma' \vdash C}{applications of (W)}
{\Theta, A \vee B, \Gamma \vdash C}}\pravilo{$\vee$L}\]
Note that this step involves permuting applications of (W) one with another.
In the sequence of applications of (W) below the sequent
$\Theta', A \vee B, \Gamma' \vdash C$ there are $max(k_1, k_2)$
applications of (W) tied to $G$ from $\Theta$, where $G$, $k_1$ and
$k_2$ are as in the formulation of the lemma.

If $n+m=0$, then this proof is W-normal.

If $n>0$, then our proof is transformed into
{\footnotesize \[\f{\f{
\Bpak{\Apak{\pi_1}{\Phi \vdash C}{\tiny applications \\ of (K)}
{\Theta', A,\ldots,A, \Gamma' \vdash C}}{\tiny $n\!-\!1$ appl. \\ of (W)}
{\Theta', A, A, \Gamma' \vdash C}
\rzbc
\Bpak{\f{\Apak{\pi_2}{\Psi \vdash C}{\tiny applications \\ of (K)}
{\Theta', B,\ldots,B, \Gamma' \vdash C}}
{\Theta', A, B,\ldots,B, \Gamma' \vdash C}\pravilo{K}}{\tiny $m$ appl. \\ of (W)}
{\Theta', A, B, \Gamma' \vdash C}
}
{\Theta', A, A \vee B, \Delta' \vdash C}\pravilo{$\vee$L}
\rzbc
\Bpak{\f{\Apak{\pi_2}{\Psi \vdash C}{\tiny applications \\ of (K)}
{\Theta', B,\ldots,B, \Gamma' \vdash C}}
{\Theta', B,\ldots,B, A \vee B, \Gamma' \vdash C}\pravilo{K}}
{\tiny $m$ appl. \\ of (W)}
{\Theta', B, A \vee B, \Gamma' \vdash C}
}
{\Dpak{\f{\Theta', A \vee B, A \vee B, \Gamma' \vdash C}
{\Theta', A \vee B, \Gamma' \vdash C}\pravilo{W}}{applications of (W)}
{\Theta, A \vee B, \Gamma \vdash C}}\pravilo{$\vee$L} \rzb \]}
Consider the subproof whose endsequent is
$\Theta', A, A \vee B, \Gamma' \vdash C$. Its measure is $n_2 + m$, where
$n_2$ is the number of applications of (W) tied to the right $A$ in the left
premise of the last rule of this subproof. The number of applications of (W)
tied to the left $A$ of the same sequent is $n_1$ and we have $n_1 + n_2 =n-1$.
We apply the induction  hypothesis to this subproof, and therefore our proof
is transformed into
\[\f{\Apak{\pi}{\Xi \vdash C}{applications \\ of (W)}
{\Theta', A, A \vee B, \Gamma' \vdash C}
\rzc
\Bpak{\f{\Apak{\pi_2}{\Psi \vdash C}{applications \\ of (K)}
{\Theta', B,\ldots,B, \Gamma' \vdash C}}
{\Theta', B,\ldots,B, A \vee B, \Gamma' \vdash C}\pravilo{K}}
{$m$ applications \\ of (W)}
{\Theta', B, A \vee B, \Gamma' \vdash C}
}
{\f{\Theta', A \vee B, A \vee B, \Gamma' \vdash C}
{\Dpak{\Theta', A \vee B, \Gamma' \vdash C}{applications of (W)}
{\Theta, A \vee B, \Gamma \vdash C}}\pravilo{W}}\pravilo{$\vee$L}\]
where in the subproofs whose endsequent is
$\Theta', A, A \vee B, \Gamma' \vdash C$ there are only $n_1$
applications of (W) tied to $A$ and there are no applications of (W)
tied to any occurrence of a formula in $\Theta'$ and $\Gamma'$.
Consider the subproof whose endsequent is
$\Theta', A \vee B, A \vee B, \Gamma' \vdash C$. Its measure is $n_1 + m$,
and by applying the induction hypothesis, our proof is transformed
into a W-normal proof with $max(k_1, k_2)$ applications of (W) tied to
$G$ from $\Theta$ in its endsequent.

We proceed quite analogously if $m>0$.
\qed
\\[0.3cm]
\lema{5.4}
{\em A proof of the form
\[\f{
\Apak{\pi_1}{\Delta' \vdash A}{applications \\ of (W)}{\Delta \vdash A}
\rzc
\Apak{\pi_2}{\Phi \vdash C}{applications \\ of (W)}{\Theta, B, \Gamma \vdash C}
}
{\Theta, \Delta, A \str B, \Gamma \vdash C}\pravilo{$\str$L}\]
where $\pi_1$ and $\pi_2$ are tailless, can be transformed into a
W-normal proof, of the same degree, of the form
\[\Apak{\pi}{\Xi \vdash C}{applications \\ of (W)}{\Theta, \Delta,
A \str B, \Gamma \vdash C}\]
where $\pi$ is tailless, and for every occurrence $G$ of a formula
in $\Theta$, if above the right premise of} $\str$L {\em in 
the former figure
there are $k$ applications of} (W) {\em tied to $G$, then 
in the latter figure
there are $k$ applications of} (W) {\em tied to $G$. The same holds for
occurrences of formulae in $\Gamma$.}
\\[0.2cm]
\dkz
Let $n$ be the number of applications of (W) tied to $B$ in the right
premise of $\str$L in the figure of the initial proof, and let the
total number of applications of (W) above this premise be $l$. We
prove the lemma by induction on $n$. Our proof is first transformed into
\[\f{\fp{\pi_1}{\Delta' \vdash A} \rzc
\f{\Cpak{\fp{\pi_2}{\Phi \vdash C}}{$n\!-\!1$ applications of \\
(W) tied to $B$}{\Theta', B, B, \Gamma' \vdash C}}
{\Theta', B, \Gamma' \vdash C}\pravilo{W}}
{\Dpak{\Theta', \Delta', A \str B, \Gamma' \vdash C}
{applications of (W) including $l\!-\!n$ applications of (W) tied to
formulae in $\Theta$ and $\Gamma$}{\Theta, \Delta, A \str B, \Gamma \vdash C}}
\pravilo{$\str$L}\]
Note that this step involves permuting applications of (W) one with another.
The transformed proof is next transformed into
\[\f{\fp{\pi_1}{\Delta' \vdash A}
\rzb
\f{\fp{\pi_1}{\Delta' \vdash A} \rzc
\Cpak{\fp{\pi_2}{\Phi \vdash C}}{$n\!-\!1$ applications \\ of (W)}
{\Theta', B, B, \Gamma' \vdash C}}
{\Theta', \Delta', A \str B, B, \Gamma' \vdash C}\pravilo{$\str$L}
}
{\Dpak{\Theta', \Delta', A \str B, \Delta', A \str B, \Gamma' \vdash C}
{applications of (C) and (W) including $l\!-\!n$ applications of (W) tied
to formulae in $\Theta$ and $\Gamma$}
{\Theta, \Delta, A \str B, \Gamma \vdash C}}\pravilo{$\str$L} \rzb \]
By the induction hypothesis, there is a W-normal proof
\[\Dpak{\fp{\pi}{\Psi \vdash C}}{applications of (W)}{\Theta', \Delta',
A \str B, B, \Gamma' \vdash C}\]
where $\pi$ is tailless, and where there are $m \leq n-1$ applications
of (W) tied to $B$ in endsequent, and no application of (W) tied to
formulae in $\Theta'$ and $\Gamma'$.

We apply again the induction hypothesis to
\[\f{\fp{\pi_1}{\Delta' \vdash A} \rzc
\Dpak{\fp{\pi}{\Psi \vdash C}}{applications of (W)}{\Theta', \Delta',
A \str B, B, \Gamma' \vdash C}}
{\Theta', \Delta', A \str B, \Delta', A \str B, \Gamma' \vdash C}
\pravilo{$\str$L}\]
and we use Lemma 5.1 to obtain a W-normal proof of
$\Theta, \Delta, A \str B, \Gamma \vdash C$. In the final proof there are
still only $l-n$ applications of (W) tied to formulae in $\Theta$
and $\Gamma$.
\qed

We can now prove the theorem that covers the first
phase of our cut-elimination procedure.
\\[0.3cm]
\teo{5.5}
{\em Every proof of a sequent in \G\ can be reduced to a W-normal proof
of the same degree of the same sequent.}
\\[0.2cm]
\dkz
We proceed by induction on the length of the proof of our sequent in \G.

If our proof is just an axiom, then this proof is W-normal.

If our sequent is proved by the following proof
\[\f{\fp{\pi}{\Delta, A, B, \Gamma \vdash C}}{\Delta, B, A, \Gamma \vdash C}
\pravilo{C}\]
then, by the induction hypothesis, there is a W-normal proof
\[\Dpak{\fp{\pi'}{\Lambda \vdash C}}{applications of (W)}
{\Delta, A, B, \Gamma \vdash C}\]
where $\pi'$ is tailless. Then we apply Lemma 5.1.

If our sequent is proved by the following proof
\[\f{\fp{\pi}{\Theta, A, A, \Gamma \vdash C}}
{\Theta, A, \Gamma \vdash C}\pravilo{W}\]
then, by the induction hypothesis, there is a W-normal proof
\[\f{\fp{\pi'}{\Theta, A, A, \Gamma \vdash C}}
{\Theta, A, \Gamma \vdash C}\pravilo{W}\]

If our sequent is proved by the following proof
\[\f{\fp{\pi}{\Theta, \Gamma \vdash C}}
{\Theta, A, \Gamma \vdash C}\pravilo{K}\]
then, by the induction hypothesis, there is a W-normal proof
\[\Dpak{\fp{\pi'}{\Lambda \vdash C}}{applications of (W)}{\Theta, \Gamma
\vdash C}\]
where $\pi'$ is tailless. Applications 
of (W) can be permuted with K as follows:
\[
\f{\f{\Theta, \Gamma \vdash C}{\Theta', \Gamma' \vdash C}
\pravilo{W}}{\Theta', A, \Gamma' \vdash C}\pravilo{K}
\rzc
\f{\f{\Theta, \Gamma \vdash C}{\Theta, A, \Gamma \vdash C}
\pravilo{K}}{\Theta', A, \Gamma' \vdash C}\pravilo{W}\]
And, by  induction on the number of applications of (W) below $\pi'$, we
prove that there is a W-normal proof of $\Theta, A, \Gamma \vdash C$.

If our sequent is proved by the following proof
\[\f{\fp{\pi}{\Delta \vdash A} \rzc \fp{\rho}{\Theta, A, \Gamma \vdash C}}
{\Theta, \Delta, \Gamma \vdash C}\pravilo{cut}\]
then, by the induction hypothesis, we have a proof
\[\f{\Apak{\pi'}{\Lambda \vdash A}{applications \\ of (W)}{\Delta \vdash A}
\rzc \Apak{\rho'}{\Phi \vdash C}{applications \\ of (W)}
{\Theta, A, \Gamma \vdash C}}{\Theta, \Delta, \Gamma \vdash C}\pravilo{cut}\]
where $\pi'$ and $\rho'$ are tailless. We push below cut
all the applications of (W) below $\pi'$ so as to obtain
\[\f{\fp{\pi'}{\Lambda \vdash A} \rzc \Apak{\rho'}{\Phi \vdash C}
{applications \\ of (W)}{\Theta, A, \Gamma \vdash C}}
{\Cpak{\Theta, \Lambda, \Gamma \vdash C}{applications \\ of (W)}
{\Theta, \Delta, \Gamma \vdash C}}\pravilo{cut}\]
Then this proof belongs to ${\cal C}(\pi', \rho')$, and we can apply Lemma 5.2.

If our sequent is proved by the following proof
\[\f{\fp{\pi}{\Theta, A, \Gamma \vdash C}}{\Theta, A \I B, \Gamma \vdash C}
\pravilo{$\I$L}\]
then, by the induction hypothesis, there is a W-normal proof
\[\Apak{\pi}{\Lambda \vdash C}{applications \\ of (W)}
{\Theta, A, \Gamma \vdash C}\]
where $\pi$ is tailless. Applications of (W) can be permuted with ($\I$L) as
follows
\[\f{\f{\Theta, A, \Gamma \vdash C}{\Theta', A, \Gamma' \vdash C}\pravilo{W}}
{\Theta', A \I B, \Gamma' \vdash C}\pravilo{$\I$L} \rzc
\f{\f{\Theta, A, \Gamma \vdash C}{\Theta, A \I B, \Gamma \vdash C}
\pravilo{$\I$L}}{\Theta', A \I B, \Gamma' \vdash C}\pravilo{W}\]
\[\f{\f{\Theta, A, A, \Gamma \vdash C}{\Theta, A, \Gamma \vdash C}
\pravilo{W}}{\Theta, A\I B, \Gamma \vdash C}\pravilo{$\I$L} \rzc
\f{\f{\Theta, A, A, \Gamma \vdash C}{\Theta, A \I B, A, \Gamma \vdash C}
\pravilo{$\I$L}}{\f{\Theta, A \I B, A \I B, \Gamma \vdash C}
{\Theta, A \I B, \Gamma \vdash C}\pravilo{W}}\pravilo{$\I$L}\]
And, by an induction analogous to that in the proof of Lemma 5.1, we show
that there is a W-normal proof of $\Theta, A \I B, \Gamma \vdash C$.
We proceed analogously for the other ($\I$L) rule, involving $B$.

If our sequent is proved by the following proof
\[\f{\fp{\pi}{\Gamma \vdash A} \rzc \fp{\rho}{\Gamma \vdash B}}
{\Gamma \vdash A \I B}\pravilo{$\I$R}\]
then, by the induction hypothesis, there are W-normal proofs
\[\Apak{\pi'}{\Gamma' \vdash A}{applications \\ of (W)}{\Gamma \vdash A}
\rzc \Apak{\rho'}{\Gamma'' \vdash B}{applications \\ of (W)}
{\Gamma \vdash B}\]
where $\pi'$ and $\rho'$ are tailless. Then we have the W-normal
proof
\[\f{
\Apak{\pi'}{\Gamma' \vdash A}{applications \\ of (K)}{\Gamma''' \vdash A}
\rzc \Apak{\rho'}{\Gamma'' \vdash B}{applications \\ of (K)}
{\Gamma''' \vdash B}}
{\Dpak{\Gamma''' \vdash A \I B}{applications of (W)}{\Gamma \vdash A \I B}}
\pravilo{$\I$R}\]

If our sequent is proved by the following proof
\[\f{\fp{\pi}{\Theta, A, \Gamma \vdash C} \rzc \fp{\rho}
{\Theta, B, \Gamma \vdash C}}
{\Theta, A \vee B, \Gamma \vdash C}\pravilo{$\vee$L}\]
we apply the induction hypothesis to $\pi$ and $\rho$, and next we apply
Lemma 5.3.

If our sequent is proved by the following  proof
\[\f{\fp{\pi}{\Gamma \vdash A}}{\Gamma \vdash A \vee B}\pravilo{$\vee$R}\]
we apply the induction hypothesis to $\pi$, and we push applications of (W)
below $\vee$R as follows
\[\f{\f{\Theta \vdash A}{\Theta' \vdash A}\pravilo{W}}
{\Theta' \vdash A \vee B}\pravilo{$\vee$R} \rzc
\f{\f{\Theta \vdash A}{\Theta \vdash A \vee B}\pravilo{$\vee$R}}
{\Theta' \vdash A \vee B}\pravilo{W}\]
Of course, we proceed analogously with the other ($\vee$R) rule, involving $B$.

If our sequent is proved by the following proof
\[\f{\fp{\pi}{\Delta \vdash A} \rzc \fp{\rho}{\Theta, B, \Gamma \vdash D}}
{\Theta, \Delta, A \str B, \Gamma \vdash D}\pravilo{$\str$L}\]
we apply the induction hypothesis to $\pi$ and $\rho$, and next we apply
Lemma 5.4.

If our sequent is proved by the following proof
\[\f{\fp{\pi}{A, \Gamma \vdash B}}{\Gamma \vdash A \str B}\pravilo{$\str$R}\]
we apply the induction hypothesis to $\pi$ to obtain the W-normal proof
\[\Dpak{\fp{\pi'}{\Phi \vdash B}}{applications of (W)}{A, \Gamma \vdash B}\]
We next push below $\str$R each of the applications of (W) not tied to $A$ in
$A, \Gamma \vdash B$.
\qed

\section{Maximal cuts}

\noindent A cut in a proof of \G\ will be called {\em maximal}
iff its rank is 2 and 
none of its premises is an axiom. A proof will be called
{\em maximalized} iff all cuts in it are maximal. We can
prove the following theorem, which covers the second phase of our
cut-elimination procedure.
\\[0.3cm]
\teo{6.1}
{\em Every W-normal proof of a sequent in \G\ can be reduced to a
maximalized W-normal proof, of the same or of a 
lower degree, of the same sequent.}
\\[0.2cm]
\dkz
It is enough to consider a W-normal proof of the form
\[ \f{\fp{\pi}{\Delta \vdash A} \rzc \fp{\rho}{\Theta, A, \Gamma \vdash C}}
{\Theta, \Delta, \Gamma \vdash C}\pravilo{cut}\]
where the cut noted in this figure is not maximal and all cuts 
in $\pi$ and $\rho$ are maximal. The rank of such  a proof
is the rank of the nonmaximal cut. We show by induction on rank that this
proof can be reduced to a maximalized W-normal proof of the same degree
of $\Theta, \Delta, \Gamma \vdash C$.

Suppose the rank of our nonmaximal cut is 2. This means that one of its
premises is an axiom. Then we eliminate this cut by standard reduction
steps, as those in (1.1) - (1.3) of the  proof of Theorem 4.1. At this point
the degree of the proof may decrease.

Suppose now the rank of our nonmaximal cut is greater than 2. In order
to decrease the rank of the proof we introduce a number of reduction steps
that decrease first the left rank. When this rank is 1, we introduce other
reduction steps that decrease the right rank. (This is opposite to Gentzen's
procedure, where the right rank is first reduced to 1. However, the matter
is not essential, and we could proceed as Gentzen did. Gentzen need
not have reduced rank to 1 on one side, before reducing the rank on the
other side -- he could as well have worked in a zig-zag manner, passing from
one side to another before reaching 1. However, for us it is essential
that the rank on one side has fallen to 1 before we attack the rank on the
other side.)

Suppose now the left rank of the nonmaximal cut above is greater than 1.
Then in addition to the standard reduction steps like those considered in (2)
of the proof of Theorem 4.1 we have the following additional reduction steps.
$$\f{\f{\fp{\pi_1}{\Delta_2 \vdash B} \rzb \fp{\pi_2}
{\Delta_1, C, \Delta_3 \vdash A}}
{\Delta_1, \Delta_2, B \str C, \Delta_3 \vdash A}
\pravilo{$\str$L} \rzc \fp{\rho}{\Theta, A, \Gamma \vdash C}}
{\Theta, \Delta_1, \Delta_2, B \str C, \Delta_3, \Gamma \vdash C}
\pravilo{cut} \leqno{(2.3)}$$
is reduced to
\[\f{\fp{\pi_1}{\Delta_2 \vdash B}  \rzc
\f{\fp{\pi_2}{\Delta_1, C, \Delta_3 \vdash A} \rzb
\fp{\rho}{\Theta, A, \Gamma \vdash C}}
{\Theta, \Delta_1, C, \Delta_3, \Gamma \vdash C}
\pravilo{cut}}
{\Theta, \Delta_1, \Delta_2, B \str C, \Delta_3, \Gamma \vdash C}
\pravilo{$\str$L} \]
The cut in the lower figure has lower rank and we may apply the induction
hypothesis to it.
$$\f{\f{\fp{\pi_1}{\Delta_2 \vdash B} \rzb \fp{\pi_2}
{\Delta_1, B, \Delta_3 \vdash A}}{\Delta_1, \Delta_2, \Delta_3 \vdash A}
\pravilo{cut} \rzc \fp{\rho}{\Theta, A, \Gamma \vdash C}}
{\Theta, \Delta_1, \Delta_2, \Delta_3, \Gamma \vdash C}
\pravilo{cut} \leqno{(2.4)}$$
is reduced to
\[\f{\fp{\pi_1}{\Delta_2 \vdash B} \rzc \f{
\fp{\pi_2}{\Delta_1, B, \Delta_3 \vdash A} \rzb
\fp{\rho}{\Theta, A, \Gamma \vdash C}}
{\Theta, \Delta_1, B, \Delta_3, \Gamma \vdash C}
\pravilo{cut}}
{\Theta, \Delta_1, \Delta_2, \Delta_3, \Gamma \vdash C}
\pravilo{cut}\]
By the induction hypothesis, the subproof of the reduced proof
ending with the right premise of the lower cut can be reduced to a
maximalized W-normal proof, of the same or of a lower 
degree, of the same sequent. The
first step of this reduction, which is one of the reduction steps
(2.1)-(2.3), makes the lower cut maximal, and subsequent steps leave it so.
We must apply (2.1)-(2.3) because the left rank of the upper cut in the
lower figure is greater than 1 (the proof $\pi_2$ cannot be an axiom, and the
right rank of the upper cut in the first figure is 1), and, moreover, $\pi_2$
cannot end with a cut. Note that in the reduction step (2.1) the rule R
cannot be (W). 

Suppose now the left rank of our cut is 1 and the
right rank is greater than 1. Then in addition to the standard reduction steps
like those considered in (3) of the proof of Theorem 4.1 (except for (3.1)
with R being (W), and (3.4), which we don't have because of
W-normality), we have the following additional cases.
$$\f{\fp{\pi}{\Delta \vdash A} \rzc \f{
\fp{\rho_1}{\Theta_2, A, \Gamma_1 \vdash B} \rzb
\fp{\rho_2}{\Theta_1, D, \Gamma_2 \vdash C}}
{\Theta_1, \Theta_2, A, \Gamma_1, B \str D, \Gamma_2 \vdash C}
\pravilo{$\str$L}}
{\Theta_1, \Theta_2, \Delta, \Gamma_1, B \str D, \Gamma_2 \vdash C}
\pravilo{cut} \leqno{(3.5)}$$
is reduced to
\[\f{ \f{\fp{\pi}{\Delta \vdash A} \rzb
\fp{\rho_1}{\Theta_2, A, \Gamma_1 \vdash B}}
{\Theta_2, \Delta, \Gamma_1 \vdash B}
\pravilo{cut} \rzc \fp{\rho_2}{\Theta_1, D, \Gamma_2 \vdash C}}
{\Theta_1, \Theta_2, \Delta, \Gamma_1, B \str D, \Gamma_2 \vdash C}
\pravilo{$\str$L} \]
We have analogous reduction steps when $A$ in the initial proof is in
$\Theta_1$ or $\Gamma_2$.
$$\f{\fp{\pi}{\Delta \vdash A} \rzc
\f{\Cpak{\fp{\rho}{C_1,\ldots,C_1,\Gamma_1, A, \Gamma_2 \vdash C_2}}
{applications \\ of (W)}{C_1, \Gamma_1, A, \Gamma_2 \vdash C_2}}
{\Gamma_1, A, \Gamma_2 \vdash C_1 \str C_2}\pravilo{$\str$R}}
{\Gamma_1, \Delta, \Gamma_2 \vdash C_1 \str C_2}\pravilo{cut} \leqno{(3.6)}$$
provided $\rho$ is tailless, is reduced to
\[\f{\fp{\pi}{\Delta \vdash A} \rzb
\fp{\rho}{C_1,\ldots,C_1,\Gamma_1, A, \Gamma_2 \vdash C_2}}
{\Cpak{C_1,\ldots,C_1,\Gamma_1, \Delta, \Gamma_2 \vdash C_2}
{applications \\ of (W)}{\f{C_1, \Gamma_1, \Delta, \Gamma_2 \vdash C_2}
{\Gamma_1, \Delta, \Gamma_2 \vdash C_1 \str C_2}
\pravilo{$\str$R}}}\pravilo{cut}\]
$$\f{\fp{\pi}{\Delta \vdash A} \rzc \f{
\fp{\rho_1}{\Theta_2, A, \Gamma_1 \vdash B} \rzb
\fp{\rho_2}{\Theta_1, B, \Gamma_2 \vdash C}}
{\Theta_1, \Theta_2, A, \Gamma_1, \Gamma_2 \vdash C}
\pravilo{cut}}
{\Theta_1, \Theta_2, \Delta, \Gamma_1, \Gamma_2 \vdash C}
\pravilo{cut} \leqno{(3.7)}$$
is reduced to
\[\f{ \f{\fp{\pi}{\Delta \vdash A} \rzb
\fp{\rho_1}{\Theta_2, A, \Gamma_1 \vdash B}}
{\Theta_2, \Delta, \Gamma_1 \vdash B}
\pravilo{cut} \rzc \fp{\rho_2}{\Theta_1, B, \Gamma_2 \vdash C}}
{\Theta_1, \Theta_2, \Delta, \Gamma_1, \Gamma_2 \vdash C}
\pravilo{cut} \]
By the induction hypothesis, the subproof of the reduced proof
ending with the left premise of the lower cut can be reduced to a
maximalized W-normal proof, of the same or of a lower degree, 
of the same sequent. The
first step of this reduction, which is one of the reduction steps
(3.1)-(3.6), makes the lower cut maximal, and subsequent steps leave it so.
We must apply (3.1)-(3.6) because the left rank of the upper cut in the
lower figure is equal to 1 and the right rank is greater than 1
(the proof $\rho_1$ cannot be an axiom, and the
left rank of the upper cut in the first figure is 1), and, moreover, $\rho_1$
cannot end with a cut.
$$\f{\fp{\pi}{\Delta \vdash A} \rzc \f{
\fp{\rho_1}{\Theta_3 \vdash B} \rzb
\fp{\rho_2}{\Theta_1, A, \Theta_2, B, \Gamma \vdash C}}
{\Theta_1, A, \Theta_2, \Theta_3, \Gamma \vdash C}
\pravilo{cut}}
{\Theta_1, \Delta, \Theta_2, \Theta_3, \Gamma \vdash C}
\pravilo{cut} \leqno{(3.8)}$$
is reduced to
\[\f{\fp{\rho_1}{\Theta_3 \vdash B} \rzc \f{
\fp{\pi}{\Delta \vdash A} \rzb
\fp{\rho_2}{\Theta_1, A, \Theta_2, B, \Gamma \vdash C}}
{\Theta_1, \Delta, \Theta_2, B, \Gamma \vdash C}
\pravilo{cut}}
{\Theta_1, \Delta, \Theta_2, \Theta_3, \Gamma \vdash C}
\pravilo{cut} \]
and we reason as for (3.7), by applying the
induction hypothesis to the subproof of the reduced
proof ending with the right premise of the lower
cut. We have an analogous reduction step
when $A$ in the initial proof is in $\Gamma$.
\qed
In terms of categories, the reduction steps (2.4) and (3.7) in the proof
above correspond to associativity of composition, whereas (3.8) corresponds
to bifunctoriality equalities.

We can now finally go into the third phase of our cut-elimination
procedure, which
is covered by the following theorem.
\\[0.3cm]
\teo{6.2}
{\em Every maximalized proof of degree greater than} 0 {\em of a sequent of
\G\ can
be reduced to a proof of lower degree of the same sequent.}
\\[0.2cm]
\dkz
Take a maximalized proof of \G\ of degree greater than 0, and starting from
the top of the proof apply to every maximal cut of the initial proof either
the standard reduction steps like those of (1.5) and (1.6) of the proof
of Theorem 4.1, or the standard 
reduction step that consists in replacing
\[\f{\f{A, \Gamma \vdash B}{\Gamma \vdash A \str B}\pravilo{$\str$R} \rzc
\f{\Delta \vdash A \rzb \Theta, B, \Xi \vdash C}
{\Theta, \Delta, A \str B, \Xi \vdash C}\pravilo{$\str$L}}
{\Theta, \Delta, \Gamma, \Xi \vdash C}\pravilo{cut}\]
by
\[\f{\f{\Delta \vdash A \rzb A, \Gamma \vdash B}{\Delta, \Gamma \vdash B}
\pravilo{cut} \rzc \afrac{\Theta, B, \Xi \vdash C}}
{\Theta, \Delta, \Gamma, \Xi \vdash C}\pravilo{cut}\]
The result is a proof whose degree  has decreased.
\qed

By applying successively the first phase, the second phase and the third
phase of our procedure, i.e. Theorems 5.5, 6.1 and 6.2, and
then again the first phase , the second phase etc., we must obtain
after one second or third phase a proof of degree 0. If this
phase was a second phase, then there are no cuts
in this proof, whereas if this phase was a third phase, then there are
cuts in the proof and all of them have atomic cut formulae.
By applying in the latter case once more the first and second phase
we will end up with a cut-free proof, because there are no
maximal cuts of degree 0.

\section{Concluding comments}

\noindent It is instructive to compare Gentzen's cut-elimination
procedure with ours at the place where Gentzen has {\em critical
mixes} (see the end of Section 2). These critical mixes correspond to
the maximal cuts whose reduction we postpone until the third phase
of our procedure. Gentzen's separation of a critical mix out of a mix,
and leaving it below, corresponds to something achieved in the first and
second phase of our procedure. When in the first phase a cut is pushed above
a contraction and is replaced by two cuts, the second phase will ensure
that the maximal cut that corresponds to the critical mix will be
at its proper place below other cuts.

To work in the presence of the lattice connectives $\I$ and $\vee$
our procedure presupposes the presence of thinning (see the proofs
of Lemma 5.3 and Theorem 5.5, case with ($\I$R)). So this procedure
as it is formulated here cannot be transferred to relevant logic,
which has contraction but lacks thinning, except if in this logic we omit
the ``additive'', i.e. lattice, connectives and restrict
ourselves to ``multiplicative'' connectives.

The problem with the lattice connectives $\I$ and $\vee$ is that in the rules
($\I$R) and ($\vee$L) there are implicit contractions: in terms of
a multiplicative rule, ($\I$R) could be reconstructed as
\[\f{\Gamma \vdash A \rzc \Gamma \vdash B}
{\Dpak{\Gamma, \Gamma \vdash A \I B}{applications of (C) and (W)}
{\Gamma \vdash A \I B}}\pravilo{multiplicative rule}\]
while for ($\vee$L) there is no such simple reconstruction, but similar
contractions are involved.
The W-normalization of the first phase of our procedure does not take
care of these implicit contractions; i.e., these are not pushed below
other rules as far as they can go. Because of that we can say that when in
the second phase of that procedure, in cases (2.2) and (3.2) of the proof
of Theorem 6.1 (which are taken over from the proof of Theorem 4.1),
there is an increase in size in the transformed proof, this
increase is again due to contraction. Contraction is, of course,
to blame for the increase in size that occurs in the first phase of the
procedure.

All the steps of our cut-elimination procedure are covered by
equalities of bicartesian closed categories, which is not the case for all the
steps of Gentzen's procedure. The categorially unjustified steps of
[Gentzen 1935] are like the following step, licenced by 3.121.1 in which
\[\f{\f{A \vdash A \rzb A \vdash A}{A, A \str A \vdash A}\pravilo{$\str$L}
\rzc \f{A \vdash A \rzb A \vdash A}{A, A \str A \vdash A}
\pravilo{$\str$L}}
{A, A \str A, A \str A \vdash A}\pravilo{mix, i.e. cut}\]
is replaced by
\[\f{\f{A \vdash A \rzb A \vdash A}{A, A \str A \vdash A}
\pravilo{$\str$L}}
{A, A \str A, A \str A \vdash A}
\pravilo{\footnotesize thinning and interchange}\]
Another problem  is that Gentzen's mix
\[\f{\Gamma \vdash A \rzb \Delta \vdash C}{\Gamma, \Delta^\ast \vdash C}\]
is strict in the sense that in $\Delta^\ast$ we must omit {\em all}
the occurrences of $A$, whereas a ``liberal'' mix where in $\Delta^\ast$
we must omit {\em some}, but not necessarily all, occurrences of $A$ is
better justified categorially. In terms of Gentzen's strict mix
the following cut
\[\f{\Gamma \vdash A \rzb A, A \vdash C}{\Gamma, A \vdash C}\pravilo{cut}\]
is reconstructed as
\[\f{\Gamma \vdash A \rzb A, A \vdash C}{\Dpak{\Gamma \vdash C}
{thinning and interchanges}{\Gamma, A \vdash C}}\pravilo{mix}\]
which is not always justified. However, it is possible to mend
Gentzen's mix-elimination procedure so that all of its steps are justified by
equalities of bicartesian closed categories.
\vspace{0.5cm}
\begin{center}
{\bf References}
\end{center}
\vskip 10pt

{\baselineskip=0.8\baselineskip
\noindent\hangindent=\parindent 
Borisavljevi\' c, M. [1997] A cut-elimination proof in
intuitionistic predicate logic. {\em Annals of Pure
and Applied Logic} {\bf 99} 105-136.

\vskip 5pt

\noindent\hangindent=\parindent 
Carbone, A. [1997] Interpolants, cut elimination and
flow graphs for the propositional calculus. {\em Annals of Pure
and Applied Logic} {\bf 83} 249-299.

\vskip 5pt

\noindent\hangindent=\parindent 
Curry, H. B. [1963]
{\em Foundations of Mathematical Logic}, McGraw Hill.

\vskip 5pt

\noindent\hangindent=\parindent 
Do\v sen, K., and Petri\' c, Z. [1999] Cartesian isomorphisms are
symmetric monoidal: A justification of linear logic. 
{\em The Journal of Symbolic Logic} {\bf 64} 227-242.

\vskip 5pt

\noindent\hangindent=\parindent 
Dyckhoff, R., and Pinto, L. [1997] Permutability of proofs in
intuitionistic sequent calculi. University of St Andrews Research Report
CS/97/7 (expanded version of a paper in
{\em Theoretical Computer Science} {\bf 212}, 1999, 141-155).

\vskip 5pt

\noindent\hangindent=\parindent 
Gentzen, G. [1935] Untersuchungen  \" uber das logische Schlie\ss en.
{\em Mathematische Zeitschrift} {\bf 39} 176-210, 405-431 (English
translation in [Gentzen 1969]).

\vskip 5pt

\noindent\hangindent=\parindent 
Gentzen, G. [1938] Neue Fassung des Wiederspuchsfreiheitsbeweises  
f\" ur die reine Zahlentheorie. {\em Forschungen zur Logik und zur
Grundlegung der exakten Wissenschaften, N.S.} {\bf 4} 19-44 
(English translation in [Gentzen 1969]).

\vskip 5pt

\noindent\hangindent=\parindent 
Gentzen, G. [1969] {\em The Collected Papers of Gerhard Gentzen},
Szabo, M.E. (ed.), North-Holland.

\vskip 5pt

\noindent\hangindent=\parindent 
Girard, J.-Y., Scedrov, A., and Scott, P. J. [1992] Bounded
linear logic: A modular approach to polynomial-time 
computability. {\em Theoretical Computer Science} {\bf 97} 1-66.

\vskip 5pt

\noindent\hangindent=\parindent 
Kleene, S.C. [1952] Permutability of inferences in Gentzen's 
calculi LK and LJ. In: Kleene, S.C. {\em Two Papers on the Predicate
Calculus}, American Mathematical Society 1-26. 

\vskip 5pt

\noindent\hangindent=\parindent 
Lambek, J. [1958] The mathematics of sentence structure. 
{\em The American Mathematical Monthly} {\bf 65} 154-170
(reprinted in: Buszkowski, W. {\em et al.} (eds) 
{\em Categorial Grammar}, Benjamins, 1988, 153-172).

\vskip 5pt

\noindent\hangindent=\parindent 
Minc, G.E. [1996] Normal forms for sequent derivations.
In: Odifreddi, P. (ed.) {\em Kreiseliana: About and Around Georg 
Kreisel}, Peters 469-492.

\vskip 5pt

\noindent\hangindent=\parindent 
Szabo, M.E. [1978] {\em Algebra of Proofs}, North-Holland.

\vskip 5pt

\noindent\hangindent=\parindent 
Zucker, J. [1974] The correspondence between cut-elimination and
normalization. {\em Annals of Mathematical Logic} {\bf 7} 1-112.

\vskip 5pt}

\end{document}